\documentclass[a4paper]{article}
\usepackage{amsmath,amsthm,amssymb}

{\theoremstyle{definition}\newtheorem{definition}{Definition}[section]
\newtheorem{notation}[definition]{Notation}
\newtheorem{notations}[definition]{Notations}
\newtheorem{convention}[definition]{Convention}

\newtheorem{remark}[definition]{Remark}
\newtheorem{example}[definition]{Example}
\newtheorem{examples}[definition]{Examples}}

\newtheorem{proposition}[definition]{Proposition}
\newtheorem{lemma}[definition]{Lemma}
\newtheorem{theorem}[definition]{Theorem}
\newtheorem{theoremquotes}[definition]{\lq Theorem\rq}
\newtheorem{corollary}[definition]{Corollary}

\renewcommand{\Im}{\operatorname{Im}}
\newcommand{\PSL}{\operatorname{PSL}}
\newcommand{\Ind}{\operatorname{Ind}}
\newcommand{\si}{\sigma}
\newcommand{\recht}{\rightarrow}
\newcommand{\Ker}{\operatorname{Ker}}
\newcommand{\Out}{\operatorname{Out}}
\newcommand{\Aut}{\operatorname{Aut}}
\newcommand{\Char}{\operatorname{Char}}
\newcommand{\cZ}{\mathcal{Z}}
\newcommand{\actson}{\curvearrowright}
\newcommand{\cR}{\mathcal{R}}

\newcommand{\al}{\alpha}
\newcommand{\ox}{\overline{x}}
\newcommand{\ok}{\overline{k}}
\newcommand{\Hom}{\operatorname{Hom}}
\newcommand{\Ad}{\operatorname{Ad}}
\newcommand{\om}{\omega}
\newcommand{\cG}{\mathcal{G}}
\newcommand{\oDelta}{\overline{\Delta}}
\newcommand{\vphi}{\varphi}
\newcommand{\Stab}{\operatorname{Stab}}
\newcommand{\cD}{\mathcal{D}}
\newcommand{\GL}{\operatorname{GL}}
\newcommand{\Z}{\mathbb{Z}}
\newcommand{\cU}{\mathcal{U}}
\newcommand{\SL}{\operatorname{SL}}
\newcommand{\Q}{\mathbb{Q}}
\newcommand{\be}{\beta}
\newcommand{\cO}{\mathcal{O}}
\newcommand{\id}{\operatorname{id}}
\newcommand{\Gal}{\operatorname{Gal}}
\newcommand{\C}{\mathbb{C}}
\newcommand{\cF}{\mathcal{F}}
\newcommand{\eps}{\varepsilon}
\newcommand{\cC}{\mathcal{C}}

\newcommand{\cL}{\mathcal{L}}

\newcommand{\Inn}{\operatorname{Inn}}
\newcommand{\M}{\operatorname{M}}
\newcommand{\ot}{\otimes}
\newcommand{\N}{\mathbb{N}}
\newcommand{\B}{\operatorname{B}}
\newcommand{\lspan}{\operatorname{span}}
\newcommand{\notembed}[1]{\underset{#1}{\not\prec}}
\newcommand{\embed}[1]{\underset{#1}{\prec}}
\newcommand{\Mtil}{\widetilde{M}}

\newcommand{\cI}{\mathcal{I}}
\newcommand{\Tr}{\operatorname{Tr}}
\newcommand{\dis}{\displaystyle}
\newcommand{\Om}{\Omega}
\newcommand{\twist}[1]{\rtimes_{#1}}
\newcommand{\Sp}{\operatorname{Sp}}
\newcommand{\R}{\mathbb{R}}
\newcommand{\setmin}{-}
\newcommand{\spect}{\operatorname{sp}}
\newcommand{\F}{\mathbb{F}}
\newcommand{\cB}{\mathcal{B}}

\begin{document}
\begin{center}
{\LARGE\bf Strong rigidity of generalized Bernoulli actions and
  computations of their symmetry groups}

\bigskip

{\sc by Sorin Popa\footnote{Supported in part by NSF Grant DMS-0601082}\footnote{\parbox[t]{15cm}{Mathematics Department; University of
    California at Los Angeles, CA 90095-1555 (United States) \\
    E-mail: popa@math.ucla.edu}\vspace{1mm}} and Stefaan Vaes\footnote{\parbox[t]{15cm}{Institut de
    Math{\'e}matiques de Jussieu; Alg{\`e}bres d'Op{\'e}rateurs; 175, rue du
    Chevaleret; F--75013 Paris (France) \\ E-mail: vaes@math.jussieu.fr}\vspace{1mm}}\footnote{Department of Mathematics; K.U.Leuven; Celestijnenlaan 200B; B--3001 Leuven (Belgium)}}
\end{center}



\begin{abstract}
\noindent We study equivalence relations and
    II$_1$ factors associated with (quotients of) generalized
    Bernoulli actions of Kazhdan groups. Specific families of these
    actions are entirely classified up to isomorphism of II$_1$
    factors. This yields explicit computations of outer automorphism
    and fundamental groups. In particular, every finitely presented
    group is concretely realized as the outer automorphism group of a
    continuous family of non stably isomorphic II$_1$ factors.
\end{abstract}

\section*{Introduction}

Given a measure preserving action of a countable group $G$ on the probability space $(X,\mu)$, Murray and von Neumann introduced the \emph{group measure space
construction} (or crossed product) denoted by $L^\infty(X,\mu) \rtimes G$, which is a von Neumann algebra acting on the Hilbert space $\ell^2(G,L^2(X,\mu))$. We are interested in three types of equivalences between actions $G
\actson (X,\mu)$ and $\Gamma \actson (Y,\eta)$, defined as follows.
\begin{enumerate}
\item {\it Von Neumann equivalence:} the associated crossed product von Neumann algebras are isomorphic.
\item {\it Orbit equivalence:} there exists a measure space isomorphism $\Delta : X \recht Y$ sending orbits to orbits almost everywhere.
\item {\it Conjugacy:} there exists a measure space isomorphism $\Delta : X \recht Y$ and a group isomorphism $\delta : G \recht \Gamma$ such that
$\Delta(g \cdot x) = \delta(g) \cdot \Delta(x)$ almost everywhere.
\end{enumerate}
It is trivial that conjugacy implies orbit equivalence and due to
observations of Dye \cite{Dye1,Dye2} and Feldman-Moore \cite{FM}, orbit equivalence implies von
Neumann equivalence. An \emph{orbit equivalence rigidity theorem} deduces for certain families of group actions, conjugacy out of their mere orbit
equivalence, while a \emph{von Neumann rigidity theorem} deduces orbit equivalence or even conjugacy of the actions out of von Neumann equivalence.
The paper consists of two parts, dealing with orbit equivalence and von Neumann rigidity theorems for families of generalized Bernoulli actions and their
quotients. We rely on techniques and results of the first
author, see \cite{P0,P1,P2}.

\subsubsection*{Orbit equivalence rigidity}

Orbit equivalence rigidity theory was initiated by the pioneering work
 of Zimmer \cite{Zim1,Zim}, who proved that orbit equivalence between
essentially free ergodic actions of lattices in higher rank simple Lie
 groups entails local isomorphism of the ambient Lie groups. In
 particular, the groups $\SL(n,\Z)$ do not admit orbit equivalent
 essentially free ergodic actions for different values of
 $n$. Zimmer's proof is based on his celebrated cocycle rigidity
 theorem, which gave rise to many future developments (see
 \cite{Fur0,Fur1,Fur2,Geft1,Geft2} to quote a few of them).
Furman developed in \cite{Fur1,Fur2} a new technique and combining it
 with Zimmer's results, obtained the first superrigidity theorems: if the action of a higher rank
lattice is orbit equivalent with an arbitrary essentially free action, both actions are conjugate (modulo finite subgroups). Another class of orbit equivalence
rigidity theorems was obtained by Monod and Shalom in
 \cite{MS}. Gaboriau proved orbit equivalence rigidity results
 of a different type, e.g.\ showing that the free groups $\F_n$ do not
 admit orbit equivalent essentially free ergodic actions for different
 values of $n$, see \cite{Gab1,Gab2}.

Using operator algebra techniques, the first author proved in
\cite{P0} a cocycle superrigidity theorem for malleable actions of $w$-rigid groups\footnote{A countable group is said to be $w$-rigid if it
admits an infinite normal subgroup with the relative property
(T).}, deducing as well orbit equivalence superrigidity for these actions. This general theorem, that applies to all generalized Bernoulli
actions, is crucial in this paper.

In the first part of the paper, we deal with orbit equivalence
rigidity theorems for quotients of generalized Bernoulli actions given by
\begin{equation} \label{eq.quotientber}
G \actson (X,\mu):=\Bigl( \prod_{I} (X_0,\mu_0) \Bigr)/K
\end{equation}
and constructed from the following data: a $w$-rigid group $G$ acting on a
countable set $I$ and a compact group $K$ acting on the probability space
$(X_0,\mu_0)$. In the construction, we make act $K$ diagonally on
$\prod_I (X_0,\mu_0)$ and take the quotient. It is shown that these
quotients of generalized Bernoulli actions are orbitally rigid: if $G
\actson (X,\mu)$ is orbit equivalent with an arbitrary essentially
free action, both actions are \lq conjugate modulo finite
groups\rq. Moreover, under good conditions, the complete data of $G
\actson I$ and $K \actson (X_0,\mu_0)$ can be recovered from the
equivalence relation given by $G$-orbits yielding a complete
classification up to orbit equivalence of certain families of
actions of the form~\eqref{eq.quotientber}. In particular, if $G$ runs
through the $w$-rigid groups without finite normal subgroups and if
$K$ runs through the non-trivial compact second countable groups, the
actions
$$G \actson \Bigl(\prod_G (K,\text{Haar}) \Bigr)/K$$
are all non stably orbit equivalent. As such, every $w$-rigid group admits
this explicit continuous family of non stably orbit equivalent actions
(cf.\ \cite{Hjo,P3}).

Self orbit equivalences yield, modulo the trivial ones, the outer automorphism group of a type II$_1$ equivalence relation. The computation of the
outer automorphism group of a concrete equivalence relation, is
usually a hard problem. The first actual computations, based on
Zimmer's work \cite{Zim1}, were done by Gefter \cite{Geft1,Geft2}, who
obtained the first equivalence relations without outer
automorphisms. A systematic treatment, including many concrete
computations, has been given by Furman \cite{Fur0}. Other examples of equivalence relations without outer automorphisms have been constructed by
Monod and Shalom \cite{MS} and by Ioana, Peterson and the first author
in \cite{IPP}. In this paper, we obtain many concrete computations of outer automorphism groups, yielding the following quite easy
example of a continuous family of non stably orbit equivalent actions
without outer automorphisms. The equivalence relation given by the
orbits of
$$(\Z^n \rtimes \GL(n,\Z)) \actson \prod_{\Z^n} (X_0,\mu_0) \quad (n \geq 2)$$
has outer automorphism group given by $\Aut(X_0,\mu_0)$, which is trivial when $\mu_0$ is atomic with different weights. Moreover, the
equivalence relation remembers the probability space
$(X_0,\mu_0)$ (i.e.\ the list of weights of the atoms and the weight
 of the continuous part). The same techniques yield continuous families of non stably isomorphic type II$_1$ equivalence
relations $\cR$ with $\Out \cR$ any prescribed countable group.

\subsubsection*{Von Neumann equivalence rigidity}

The first rigidity phenomena in von Neumann algebra theory were
discovered by Connes in \cite{C1}. He showed that whenever $\Gamma$ is
a property~(T) group with infinite conjugacy classes (ICC), the group von Neumann algebra $\cL(\Gamma)$
has countable outer automorphism group and countable fundamental
group. Without being exhaustive, we cite the following von Neumann
rigidity results obtained during the last 25 years:
\cite{C1,C2,C10,CJ,CH,GN,P4,P5,P8}.

The first von Neumann strong rigidity theorem, deducing conjugacy of
actions out of their von Neumann equivalence, was recently obtained by
the first author \cite{P1,P2}, who proved that von Neumann equivalence
of an essentially free action of a $w$-rigid group on the one hand and
the Bernoulli action of an ICC group on the other hand, entails
conjugacy of the actions through an isomorphism of the groups. In
particular, if $\Z \wr G := \Z^{(G)} \rtimes G$ denotes the wreath
product, the II$_1$ factors $\cL(\Z \wr G)$ are non-isomorphic for
non-isomorphic $w$-rigid ICC groups $G$. This was the first result
embedding a large class of groups \lq injectively\rq\ into the
category of II$_1$ factors. It proves a relative version of Connes's
conjecture \cite{C5} that the isomorphism of $\cL(\Gamma)$ and
$\cL(G)$ implies (virtual) isomorphism of $\Gamma$ and $G$ whenever
$\Gamma$ is a property~(T) ICC group.

In the second part of the paper, we prove a von Neumann rigidity
theorem for generalized Bernoulli actions of the form
\begin{equation}\label{eq.genber}
G \actson \prod_{G/G_0} (X_0,\mu_0) \; ,
\end{equation}
where $G$ is a $w$-rigid ICC group and where $G_0 \subset G$ is a \lq very non-normal\rq\ subgroup. In fact, under good conditions, the crossed product von Neumann
algebra remembers all the data: the inclusion $G_0 \subset G$ and the probability space $(X_0,\mu_0)$. This allows to fully classify certain families of generalized Bernoulli actions up to von Neumann equivalence.

We also compute the outer automorphism groups of II$_1$ factors
associated with certain actions of the form~\eqref{eq.genber}. Note that the
first author's von Neumann strong rigidity theorem \cite{P2} for Bernoulli
actions $G \actson (X,\mu):=\prod_G (X_0,\mu_0)$ allows to describe
the outer automorphism group of the crossed product in terms of
the normalizer of $G \subset \Aut(X,\mu)$ (which contains the measure
space automorphisms of $(X,\mu)$ that commute with the Bernoulli
action). But the actual computation of such normalizers for Bernoulli
actions remained open. For actions of the form~\eqref{eq.genber} and
$G_0 \subset G$ sufficiently non-normal, the normalizer of $G \subset
\Aut(X,\mu)$ can be computed. This leads to continuous families
of non stably isomorphic II$_1$ factors $M$ such that $\Out M$ is any
prescribed group of finite presentation and we strongly believe that
the same method permits in fact to get any countable grou. In particular, we get concrete
continuous families of II$_1$ factors $M$ with trivial outer
automorphism group and trivial fundamental group, as follows: let $G =
\Z^4 \rtimes \SL(4,\Z)$, with subgroup $G_0$ consisting of the elements $\pm A^n, n \in
\Z$, where $A \in \SL(4,\Z)$ is given by
\begin{equation}\label{eq.matrixA}
A=\begin{pmatrix} 0 & 0 & -1 & -1 \\ 1 & 0 & 0 & 0 \\ 0 & 1 & 0 & 0
\\ 0 & 0 & 1 & 0\end{pmatrix} \; .
\end{equation}
If $G \actson (X,\mu)$ is given by~\eqref{eq.genber}, the outer
automorphism group of $M:=L^\infty(X,\mu) \rtimes G$ is given by
$\Aut(X_0,\mu_0)$, which is trivial when $\mu_0$ is atomic with
 different weights. Also, $M$ remembers $(X_0,\mu_0)$.
We mention here that the first examples
of II$_1$ factors $M$ with $\Out M$ trivial, were given by Ioana,
Peterson and the first author in \cite{IPP}. The construction in
\cite{IPP} uses the Baire category theorem, yielding an existence theorem.

Adding scalar $2$-cocycles, we also obtain continuous families of
II$_1$ factors $M$ with trivial fundamental group, trivial $\Out M$
and with $M \not\cong M^\text{op}$. The first II$_1$ factors not
anti-isomorphic to themselves were constructed by Connes in
\cite{C10}, distinguishing $M$ and $M^\text{op}$ using his
$\chi(M)$-invariant. These examples however have fundamental group
$\R^*_+$ and large outer automorphism group since they are McDuff factors.

Reinterpreting the crossed product by a generalized Bernoulli action as the group algebra of a wreath product group, we provide the first examples of
ICC groups $G$ such that
\begin{equation}\label{eq.millenium}
\Out(\cL(G)) \cong \Char(G) \rtimes \Out(G) \; ,
\end{equation}
with respect to the obvious homomorphism mapping $(\om,\delta) \in \Char(G) \rtimes \Aut(G)$ to the automorphism of $\cL(G)$ defined by
$\theta_{\om,\delta}(u_g) = \om(g) u_{\delta(g)}$. This is related to one of Jones' millenium problems: following up Connes' rigidity conjecture
\cite{C5}, Jones asked in \cite{JoMill} if the isomorphism \eqref{eq.millenium} holds for ICC groups having property~(T).

The von Neumann rigidity theorems in this paper are natural
continuations of the work done by the first author in \cite{P1,P2}, where a
von Neumann strong rigidity theorem for Bernoulli actions of $w$-rigid
groups is obtained. Such Bernoulli actions are \emph{mixing}, while the
generalized Bernoulli actions of the form~\eqref{eq.genber} are only
\emph{weakly mixing}. On the one hand, the mixing property was a
crucial ingredient in the proofs of \cite{P1,P2} and on the other
hand, it is the absence of mixing that makes it possible in the
current paper to give classification results and complete computations
of outer automorphism groups. Although it may seem as a technical
issue, the step from mixing to weakly mixing actions in the von
Neumann rigidity theorem, was our major challenge.

\subsubsection*{Organization of the paper}

In the first section, we introduce some terminology and give the
statements of our main theorems. In Section \ref{sec.generalities}, we
deal with notations and a few generalities, while in Section
\ref{sec.cocycle-superrigid} the cocycle superrigidity theorem of
\cite{P0} is recalled. Our orbit equivalence rigidity and
classification results are proved in Sections
\ref{sec.orbit-equivalence-superrigid} and
\ref{sec.orbit-equivalence-classif}. We show the von Neumann strong
rigidity theorem for generalized Bernoulli actions in Section
\ref{sec.vNstrong}. The final Section \ref{sec.examples} is devoted to
concrete computations of outer automorphism groups of equivalence
relations and II$_1$ factors.

\subsubsection*{Acknowledgment}
The second author would like to thank Bachir
Bekka for his help in constructing almost malnormal subgroups of
$\GL(n,\Z)$ and Narutaka Ozawa for the discussions about hyperbolic groups. He also thanks the University of California at Los
Angeles for their hospitality during the work on this paper.

\section{Statement of the main results}

Let $G$ be a countable group acting on the standard probability space
$(X,\mu)$ in a measure preserving way. The action $G \actson (X,\mu)$
is said to be \emph{(essentially) free} if almost every $x \in X$ has
a trivial stabilizer and is said to be \emph{ergodic} if a measurable
$G$-invariant subset has measure $0$ or $1$.

Associated with $G \actson (X,\mu)$ are the \emph{equivalence relation} $\cR$
given by $G$-orbits and the \emph{crossed product von Neumann algebra}
$M:=L^\infty(X,\mu) \rtimes G$ acting on the Hilbert space $L^2(X,\mu)
\ot \ell^2(G)$ and generated by a copy of $L^\infty(X,\mu)$ as well as a copy
of $\Gamma$ acting as follows:
$$F (\xi \ot \delta_g) = (F(g^{-1} \cdot )\xi) \ot \delta_g
\quad\text{for}\;\; F \in L^\infty(X,\mu) \quad\text{and}\quad u_h(\xi
\ot \delta_g) = \xi \ot \delta_{gh^{-1}} \quad\text{for}\;\; h \in
\Gamma \; .$$
A \emph{type II$_1$ equivalence relation} $\cR$ is an ergodic equivalence
relation on a standard non-atomic probability space $(X,\mu)$ with countable
equivalence classes preserving the probability measure $\mu$. A
\emph{type II$_1$ factor} is an infinite dimensional von Neumann
algebra with trivial center that admits a normal tracial state. If $G
\actson (X,\mu)$ is essentially free and ergodic (and if $G$ is
infinite, i.e.\ $\mu$ non-atomic), the associated
equivalence relation and the associated crossed product von Neumann
algebra are of type II$_1$.

Actions $G \actson (X,\mu)$ and $\Gamma
\actson (Y,\eta)$ are said to be \emph{orbit equivalent} if their
associated equivalence relations $\cR_G$ and $\cR_\Gamma$ are
isomorphic. They are said to be \emph{von Neumann equivalent} if the
associated crossed products are isomorphic. Note that two such actions
are orbit equivalent if and only if there exists an isomorphism $\pi :
L^\infty(X,\mu) \rtimes G \recht L^\infty(Y,\eta) \rtimes \Gamma$
sending $L^\infty(X,\mu)$ onto $L^\infty(Y,\eta)$, see \cite{FM}.

If $\cR$
is a type II$_1$ equivalence relation, the \emph{amplification} $\cR^t$ is
well defined, up to isomorphism, for all $t > 0$. If $0<t\leq 1$, one
defines $\cR^t$ as the restriction of $\cR$ to a set of measure
$t$. The \emph{fundamental group} $\cF(\cR)$ of the type II$_1$ equivalence
relation $\cR$ is defined as the group of $t > 0$ such that $\cR^t$ is
isomorphic with $\cR$. Two essentially free ergodic actions are said to be \emph{stably orbit
  equivalent} if the associated equivalence relations admit isomorphic
amplifications.

In the same way, the amplification $M^t$ of a II$_1$ factor $M$ is
defined for every $t > 0$ and the fundamental group $\cF(M)$ is
defined as the group of $t > 0$ such that $M \cong M^t$.

We denote by $\Aut(X,\mu)$ the group of measure space isomorphisms of
$(X,\mu)$ preserving the measure $\mu$, where we identify
transformations equal almost everywhere. Then, $\Aut(X,\mu)$ is
canonically isomorphic with the group of trace preserving
automorphisms of $L^\infty(X,\mu)$.

If $\cR$ is a type II$_1$ equivalence relation on $(X,\mu)$, the
group $\Aut \cR$ is defined as the group of $\Delta \in \Aut(X,\mu)$ satisfying
$(\Delta(x),\Delta(y)) \in \cR$ for almost all $(x,y) \in \cR$. The
\emph{inner automorphism group} $\Inn \cR$, also called \emph{full group} of
$\cR$, consists of $\Delta \in \Aut(X,\mu)$ such that $(x,\Delta(x))
\in \cR$ for almost all $x \in X$. The \emph{outer automorphism group}
of $\cR$ is denoted by $\Out \cR$ and defined as the quotient $\Aut
\cR / \Inn \cR$.

If $M$ is a II$_1$ factor, we denote by $\Inn M$ the group of
automorphisms of $M$ of the form $\Ad u$ for some unitary $u \in
M$. These are called inner automorphisms and the quotient $\Aut M /
\Inn M$ is called the \emph{outer automorphism group} of the II$_1$
factor $M$.

\subsection{Orbit equivalence rigidity}

\begin{definition}
Consider a countable group $G$ acting on a countable set $I$. Let $(X_0,\mu_0)$ be a non-trivial standard probability space. The \emph{$(G \actson
I)$-Bernoulli action with base $(X_0,\mu_0)$} is defined as
$$G \actson \prod_{i \in I} (X_0,\mu_0) \; .$$
We always assume that $G \actson I$ is sufficiently free: if $g \neq e$, there are infinitely many $i \in I$ with $g \cdot i \neq i$.
\end{definition}

Suppose now that $K$ is a second countable compact group acting on the base space $(X_0,\mu_0)$, preserving the probability measure $\mu_0$. Consider
$$(X,\mu) = \prod_{i \in I} (X_0,\mu_0) \qquad\text{with the diagonal action}\quad K \actson (X,\mu) \; .$$
\begin{notation}
We denote by $(X^K,\mu^K)$ the quotient of $(X,\mu)$ by the action of $K$. We view $(X^K,\mu^K)$ through the equality $(L^\infty(X^K),d\mu^K) = (L^\infty(X)^K,d\mu)$, where
the latter consists of $K$-invariant elements of $L^\infty(X)$.
\end{notation}
Let $G \actson I$ and consider the generalized Bernoulli action $G
\actson (X,\mu)$. This action obviously commutes with the action of $K$ and we consider
the action $G \actson X^K$ of $G$ on the quotient space $X^K$.

The following is an orbit superrigidity theorem for quotients of generalized Bernoulli actions. It is a corollary to Theorem \ref{thm.superrigid}
proven below.

\begin{theorem} \label{thm.orbital-superrigid}
Let $G$ be a countable group that admits an infinite normal subgroup $H$ with the relative property~(T). Assume that $G$ does not have finite normal
subgroups. Let $G \actson I$ and suppose that $H \cdot i$ is infinite for all $i \in I$. Let $K \actson (X_0,\mu_0)$ and consider the quotient of the
generalized Bernoulli action $G \actson X^K$ as above.

Any essentially free action that is orbit equivalent with $G \actson X^K$, is conjugate with $G \actson X^K$.
\end{theorem}

We now specify to actions $G \actson I$, where $I$ is the coset space $G/G_0$ for a subgroup $G_0 \subset G$. Under the right
conditions, the (quotients of) $(G \actson G/G_0)$-Bernoulli actions can be completely classified. In particular, the outer automorphism groups of
the associated equivalence relations can be computed. The \lq good\rq\ conditions on $G_0 \subset G$ are gathered in the following definition.

\begin{definition}[Condition $\cB$] \label{def.conditionB}
We say that $G_0 \subset G$ satisfies \emph{condition $\cB$} if the following holds.
\begin{itemize}
\item $G_0$ is an infinite subgroup of the countable group $G$ and $G$ does not have finite normal subgroups.
\item For all $g \in G$ with $g \not\in G_0$, the subgroup $gG_0 g^{-1} \cap G_0$ has infinite index in $G_0$. In words: the quasi-normalizer of $G_0$
in $G$ equals $G_0$.
\item Every $g \neq e$ moves infinitely many cosets $x \in G/G_0$.
\item $G$ admits an infinite normal subgroup $H$ with the relative property~(T). Moreover, the subgroup $H \cap G_0$ has infinite index in $H$.
\end{itemize}
\end{definition}

\begin{theorem} \label{thm.classification}
Denote by $\si(G_0 \subset G, K \actson (X_0,\mu_0))$ the quotient of the generalized Bernoulli action given by
$$G \actson \Bigl(\prod_{G/G_0} (X_0,\mu_0)\Bigr)^K \; .$$
Denote by $\cR(G_0 \subset G,K \actson (X_0,\mu_0))$ the equivalence relation given by the $G$-orbits. Suppose that $G_0 \subset G$ satisfies
condition $\cB$ and that $K \recht \Aut(X_0,\mu_0)$ is injective.
\begin{enumerate}
\item Any essentially free action that is orbit equivalent with $\si(G_0 \subset G,K \actson (X_0,\mu_0))$ is conjugate with the latter.
\item $\cR:=\cR(G_0 \subset G,K \actson (X_0,\mu_0))$ has trivial
  fundamental group, $\Inn \cR$ is closed in $\Aut \cR$ and
$$\Out(\cR(G_0 \subset G,K \actson (X_0,\mu_0))) =
\Hom(G/G_0 \recht \cZ(K)) \rtimes \Bigl(\frac{\Aut(G_0 \subset G)}{\Ad G_0} \times \frac{\Aut^*(K \actson
  X_0)}{K}\Bigr) \; ,$$
where $\Hom(G/G_0 \recht \cZ(K))$ denotes the group of homomorphisms from $G$ to the center of $K$ that are constantly $e$ on $G_0$ and where
$\Aut^*(K \actson X_0)$ denotes the group of measure space
automorphisms of $(X_0,\mu_0)$ that normalize $K \subset \Aut(X_0,\mu_0)$.
\item If also $\Gamma_0 \subset \Gamma$ satisfies condition $\cB$, the following are equivalent.
\begin{itemize}
\item The actions $\si(G_0 \subset G,K \actson (X_0,\mu_0))$ and $\si(\Gamma_0 \subset \Gamma,L \actson (Y_0,\eta_0))$ are stably orbit equivalent.
\item We have $(G_0 \subset G) \cong (\Gamma_0 \subset \Gamma)$ and the actions $K \actson (X_0,\mu_0)$ and $L \actson (Y_0,\eta_0)$ are conjugate through
an isomorphism $K \recht L$.
\end{itemize}
\end{enumerate}
\end{theorem}
So, the equivalence relation $\cR(G_0 \subset G,K \actson (X_0,\mu_0))$ remembers all the ingredients: the inclusion $G_0 \subset G$ and the action
$K \actson (X_0,\mu_0)$.

We get the following examples of equivalence relations with trivial outer automorphism groups (see Subsection \ref{subsec.outer-equivalence} for details).

\begin{example}
Consider $G:=\Z^n \rtimes \GL(n,\Z)$ acting on $\Z^n$ by $(x,g) \cdot
y = x + g \cdot y$ and let $\cR(X_0,\mu_0)$ be the equivalence relation given
by the orbits of
$$G \actson \Bigl(\prod_{\Z^n} (X_0,\mu_0)\Bigr) \; .$$
Then, $\Out(\cR(X_0,\mu_0)) = \Aut(X_0,\mu_0)$, which is trivial when $\mu_0$ is atomic with distinct weights. Also, $\cR(X_0,\mu_0)$ and
$\cR(Y_0,\eta_0)$ are stably isomorphic if and only if the probability spaces $(X_0,\mu_0)$ and $(Y_0,\mu_0)$ are isomorphic.
\end{example}

In Subsection \ref{subsec.examples-groups} lots of examples of $G_0 \subset G$ satisfying condition $\cB$ are constructed. In Theorem
\ref{thm.lot-of-out-equiv}, this yields, for any countable group $Q$ and any second countable compact group $K$, continuous families of type II$_1$
equivalence relations $\cR$ with $\Out(\cR) \cong Q \times \Aut(K)$.

\subsection{Strong rigidity for von Neumann algebras}

We also prove a von Neumann algebraic version of Theorem
\ref{thm.classification}, but this is much harder. We provide explicit continuous families of II$_1$ factors with trivial outer automorphism group and
trivial fundamental group. We provide quite large families of II$_1$
factors with entirely computable outer automorphism groups. In
particular, in Theorem \ref{thm.lot-of-out}, we obtain for any group
$Q$ of finite presentation, a continuous family
of II$_1$ factors $M$ with $\Out M \cong Q$. Using twisted crossed
products by a scalar $2$-cocycle, we obtain continuous families of
II$_1$ factors $M$ that are not anti-isomorphic to themselves, that
have $\Out M$ trivial and trivial fundamental group (see Theorem \ref{thm.without-anti}).

In order to prove a von Neumann algebraic version of Theorem \ref{thm.classification}, one has to strengthen condition $\cB$ (Def.\
\ref{def.conditionB}). A rude way consists in looking at the following type of groups.

\begin{definition}[Condition $\cD$] \label{def.conditionD}
We say that $G_0 \subset G$ satisfies \emph{condition $\cD$} if the following holds.
\begin{itemize}
\item $G_0$ is an infinite subgroup of the infinite conjugacy class
  group $G$.
\item $G$ admits an infinite normal subgroup with the property~(T) of
  Kazhdan (not the relative property~(T)~!).
\item For all $g \in G$ with $g \not\in G_0$, the group $g G_0 g^{-1}
  \cap G_0$ is finite.
\item $G_0$ has the Haagerup property.
\end{itemize}
\end{definition}

Note that a subgroup $G_0 \subset G$ satisfying $g G_0 g^{-1}
  \cap G_0$ finite for all $g \in G,g \not\in G_0$, is
  called \emph{almost malnormal}. Replacing finite by trivial, one
  gets the definition of a \emph{malnormal} subgroup.

The following statement is vague, but applies as such under condition
$\cD$. In Subsection \ref{subsec.outer-vN}, other large families of
groups $G_0 \subset G$ satisfying the needed conditions are
presented.

\begin{theoremquotes}
Let $G_0 \subset G$ and $\Gamma_0 \subset \Gamma$ satisfy \emph{an
  appropriate strengthening of condition $\cB$} (e.g.\ suppose that
  $(G,G_0)$ and $(\Gamma,\Gamma_0)$ satisfy condition $\cD$). Let $G \actson (X,\mu)$ be the $(G \actson
G/G_0)$-Bernoulli action with base $(X_0,\mu_0)$ and let $\Gamma \actson (Y,\eta)$ be the $(\Gamma \actson \Gamma/\Gamma_0)$-Bernoulli action with
base $(Y_0,\eta_0)$.

If $t > 0$ and
$$\pi : L^\infty(X,\mu) \rtimes G \recht \bigl( L^\infty(Y,\eta)
\rtimes \Gamma \bigr)^t$$
is a $^*$-isomorphism, then $t=1$ and there exists
\begin{itemize}
\item an isomorphism $\delta : G \recht \Gamma$ with $\delta(G_0) =
  \Gamma_0$;
\item a measure space isomorphism $\Delta_0 : (X_0,\mu_0) \recht (Y_0,\eta_0)$;
\item a character $\om \in \Char(G)$;
\item a unitary $u \in  L^\infty(Y,\eta) \rtimes \Gamma$;
\end{itemize}
such that $((\Ad u) \circ \pi)( a u_g ) = \om(g) \; \al(a) \; u_{\delta(g)}$ for all $a \in L^\infty(X,\mu)$ and $g \in G$, where $\al :
L^\infty(X,\mu) \recht L^\infty(Y,\eta)$ is given by
$$\al(a) = a \circ \Delta^{-1} \quad\text{and}\;\;
\Delta(x)_{\overline{\delta(g)}} = \Delta_0(x_{\overline{g}}) \; .$$
In particular, $M:= L^\infty(X,\mu) \rtimes G$ has trivial fundamental
group, $\Inn M$ is closed in $\Aut M$ and
$$\Out M  \cong \Bigl(\Char G \rtimes
\frac{\Aut(G_0 \subset G)}{\Ad G_0} \Bigr) \times \Aut(X_0,\mu_0) \; .$$
\end{theoremquotes}

The previous theorem can be generalized in two directions: first of all, the same result holds for generalized Bernoulli actions with non-commutative
base, i.e.\ replacing systematically $L^\infty(X_0,\mu_0)$ by the hyperfinite II$_1$ factor or by the matrix algebra $\M_n(\C)$, see Theorem
\ref{thm.vNstrong-factor}. Secondly, a version of the previous theorem holds for crossed products twisted by a scalar $2$-cocycle, see Theorems
\ref{thm.vNstrong-cocycle} and \ref{thm.without-anti}.

We obtain the following explicit examples of II$_1$ factors with trivial outer automorphism group (see Subsection \ref{subsec.outer-vN} for details and
further computations).

\begin{example}
Set $G = \Z^4 \rtimes \SL(4,\Z)$ and let $A \in \SL(4,\Z)$ be the
matrix given by \eqref{eq.matrixA} above. Let $G_0$ be the subgroup of
$G$ given by the elements $\pm A^n$, $n \in \Z$. Then, $G_0 \subset G$ satisfies condition $\cD$, $\Aut(G_0 \subset G) = \{\id\}$ and $\Char G$ is
trivial. Whenever $(X_0,\mu_0)$ is an atomic probability space, the II$_1$ factor
$$L^\infty\Bigl(\prod_{G/G_0} (X_0,\mu_0) \Bigr) \rtimes G$$
has trivial outer automorphism group and trivial fundamental
group. Moreover, they are non-stably isomorphic for different
atomic probability spaces $(X_0,\mu_0)$.
\end{example}

\subsection{Conventions}

Throughout the paper, all actions $G \actson (X,\mu)$ and all measure
space isomorphisms are assumed to
be measure preserving. All groups are supposed to be second countable.

\section{Generalized Bernoulli actions and their quotients} \label{sec.generalities}

We gather several either easy, either well known results on generalized Bernoulli actions and their quotients. Recall the following definition.

\begin{definition}
The probability measure preserving action $G \actson (X,\mu)$ is said to be \emph{weakly mixing} if
one of the following equivalent conditions is satisfied.
\begin{itemize}
\item The diagonal action $G \actson (X \times X,\mu \times \mu)$ is
  ergodic.
\item The diagonal action $G \actson (X \times Y,\mu \times \eta)$ is
  ergodic for all ergodic $G \actson (Y,\eta)$.
\item The constant functions form the only finite-dimensional
  $G$-invariant subspace of $L^2(X,\mu)$.
\item For every finite subset $\cF\subset L^\infty(X,\mu)$, there exists a
  sequence $g_n$ in $G$ such that
$$\int_X F_1(x) F_2(g_n \cdot x) \; d\mu(x) \recht \int_X F_1 \; d\mu \;\;
\int_X F_2 \; d\mu \qquad\text{for all}\;\; F_1,F_2 \in \cF \; . $$
\end{itemize}
\end{definition}

Let $G \actson (X,\mu)$ be the $(G \actson I)$-Bernoulli action with base $(X_0,\mu_0)$. If $(X_0,\mu_0)$ is non-atomic, the action $G \actson
(X,\mu)$ is essentially free if and only if every $g \neq e$ moves at least one element of $I$. If $(X_0,\mu_0)$ has atoms, the action $G \actson
(X,\mu)$ is essentially free if and only if every $g \neq e$ moves infinitely many elements of $I$. For that reason, we agree on the following
convention.

\begin{convention} \label{conv.free}
Whenever a countable group $G$ acts on the countable set $I$, we make the following freeness assumption: for every $g \neq e$, there exist infinitely
many $i \in I$ satisfying $g \cdot i \neq i$. In the same run, the group $G$ and the set $I$ are implicitly assumed to be infinite.
\end{convention}

\begin{proposition}
Let $G \actson (X,\mu)$ be the $(G \actson I)$-Bernoulli action with
base $(X_0,\mu_0)$. Then, the following conditions are equivalent.
\begin{enumerate}
\item Every orbit of $G \actson I$ is infinite.
\item $G \actson (X,\mu)$ is weakly mixing.
\item $G \actson (X,\mu)$ is ergodic.
\end{enumerate}
\end{proposition}

The non-trivial implication is 1 $\Rightarrow$ 2, which is a consequence of the following folklore lemma.

\begin{lemma} \label{lem.bete}
Let a group $G$ act on a set $I$. Then, the following conditions are equivalent.
\begin{itemize}
\item Every orbit of $G \actson I$ is infinite.
\item $G \actson I$ is weakly mixing: for every $A,B \subset
  I$ finite, there exists $g \in G$ satisfying
$$g \cdot A \cap B = \emptyset \; ,$$
\end{itemize}
\end{lemma}
\begin{proof}
Suppose that $A,B \subset I$ are finite and $g \cdot A \cap B \neq
\emptyset$ for all $g \in G$. We have to prove that $G \actson I$ admits
a finite orbit. Denote by $\chi_A$ and $\chi_B$ the indicator
functions of $A$ and $B$ and consider them as vectors in
$\ell^2(I)$. Consider the obvious unitary representation of $G$ on
$\ell^2(I)$ given by $\pi(g) \delta_i = \delta_{g \cdot i}$. Let $\xi
\in \ell^2(I)$ be the unique element of minimal norm in the closed
convex hull of $\{ \pi(g) \chi_A \mid g \in G \}$. For every $g \in
G$, we have
$$\langle \pi(g) \chi_A,\chi_B \rangle = \langle \chi_{g \cdot A} ,
\chi_B \rangle = \# (g \cdot A \cap B) \geq 1 \; .$$
It follows that $\langle \xi,\chi_B \rangle \geq 1$ and hence, $\xi
\neq 0$. By uniqueness of $\xi$, we have $\pi(g) \xi = \xi$ for all $g
\in G$, i.e.\ $\xi(g \cdot i) = \xi(i)$ for all $g \in G$, $i \in
I$. Taking $\eps > 0$ small enough,
$$I_0 := \{ i \in I \mid |\xi(i)| > \eps \}$$
is non-empty, $G$-invariant and finite because $\xi \in \ell^2(I)$.
\end{proof}

We end this section with a few observations about quotients of
generalized Bernoulli actions. So, suppose that a compact second
countable group $K$ acts on $(X_0,\mu_0)$ and let $G \actson (X,\mu)$
be the $(G \actson I)$-generalized Bernoulli action. Consider as
above the action $G \actson X^K$ on the quotient of $X$ by $K$ acting
diagonally on $(X,\mu)$.

We first observe that $G \actson X^K$ is essentially free when $G \actson I$ satisfies the freeness convention \ref{conv.free}. Indeed, let $g \neq e$ and
take $I_1 \subset I$ infinite such that $(g \cdot I_1) \cap I_1 = \emptyset$. Set $I_2 = g \cdot I_1$. We then consider the projection map
$$X^K \recht X_1^K \times X_2^K \quad\text{where}\quad X_i :=
\prod_{I_i} (X_0,\mu_0) \; .$$ Every element $x \in X^K$ fixed by $g$ is mapped to a diagonal element of the form $(y,g \cdot y)$ and the diagonal is
of measure zero because $X_1^K$ has no atoms.

The diagonal infinite product action of a compact group action satisfies an essential freeness property, which is again a folklore result. We use the
following terminology: an action $K \actson (X_0,\mu_0)$ is said to be faithful if acting by $k \neq e$ is not almost everywhere the identity
transformation. This means that the homomorphism $K \recht \Aut(X_0,\mu_0)$ is injective. If this is not the case, we can of course pass to the
quotient.

\begin{lemma} \label{lem.freeness}
Let $K$ be a compact second countable group. Let $(X_0,\mu_0)$ be a non-trivial standard probability space and $K \actson (X_0,\mu_0)$ a faithful
action. Consider
$$K \actson (X,\mu) = \prod_{i \in I} (X_0,\mu_0)$$
diagonally for some countably infinite set $I$. The action $K \actson
(X,\mu)$ is essentially free in the sense that there exists a measurable map $u : X \recht
K$ satisfying $u(k \cdot x) = k u(x)$ almost everywhere.
\end{lemma}
\begin{proof}[Sketch of proof]
It suffices to prove that there exists a subset $W \subset X$ of measure zero such that whenever $x \in X$ has a non-trivial stabilizer, then $x \in
W$. To prove this, it suffices to show that for every $k \neq e$ in $K$, there exists a neighborhood $\cU$ of $k$ in $K$ and a subset $W \subset X$
of measure zero such that every element $x \in X$ stabilized by an element of $\cU$ is contained in $W$.

Let $k \neq e$ in $K$. Since $K$ acts faithfully on $(X_0,\mu_0)$, we can find a neighborhood $\cU$ of $k$ in $K$ and a non-negligible subset $V_0
\subset X_0$ such that $(h \cdot V_0) \cap V_0 = \emptyset$ for all $h \in \cU$. But then, $W := \prod_I (X_0 \setmin V_0)$ is a set of measure zero
with the required properties.
\end{proof}

\section{The cocycle superrigidity theorem} \label{sec.cocycle-superrigid}

We recall the cocycle superrigidity theorem due to the first author \cite{P0}. We formulate it here for generalized Bernoulli actions, but it holds
more generally for \emph{malleable} actions, see \cite{P0} for details.

\begin{definition}
A Polish group $\cG$ is said to be of \emph{finite type} if $\cG$ can
be realized as a closed subgroup of the unitary group of a finite von
Neumann algebra (with separable predual).
\end{definition}

The following groups are of finite type: countable discrete groups, compact second countable groups and their direct products. These are the cases
that are needed in this paper.

Recall that, given an action $G \actson (X,\mu)$, a \emph{$1$-cocycle} with values in a Polish group $\cG$ is a measurable map
$$\om : G \times X \recht \cG \quad\text{satisfying}\quad \om(gh,x) = \om(g,h \cdot x) \; \om(h,x) \quad\text{almost everywhere.}$$
Two $1$-cocycles $\om_1$ and $\om_2$ are said to be \emph{cohomologous} if there exists a measurable map $\vphi : X \recht \cG$ satisfying
$$\om_2(g,x) = \vphi(g\cdot x) \om_1(g,x) \vphi(x)^{-1} \quad\text{almost everywhere.}$$
Note that $1$-cocycles not depending on the space variable, precisely
are homomorphisms $G \recht \cG$.

The importance of $1$-cocycles in orbit equivalence theory, stems from
the following crucial observation of Zimmer. If $\Delta : X \recht Y$
is an orbit equivalence between essentially free actions $G \actson
(X,\mu)$ and $\Gamma \actson (Y,\eta)$, the equation $\Delta(g \cdot
x) = \om(g,x) \cdot \Delta(x)$ almost everywhere, defines a
$1$-cocycle $\om$ for the action $G \actson (X,\mu)$ with values in
$\Gamma$. We call $\om$ the \emph{Zimmer $1$-cocycle} associated with
the orbit equivalence $\Delta$.

\begin{theorem}[Cocycle superrigidity theorem \cite{P0}] \label{thm.cocycle-superrigid}
Let $G$ be a countable group with an infinite normal subgroup $H$ with
the relative property~(T). Let $G \actson (X,\mu)$ be a generalized
Bernoulli action and suppose that its restriction to $H$ is weakly
mixing.

Every $1$-cocycle for $G \actson (X,\mu)$ with values in a Polish
group of finite type $\cG$ is cohomologous to a homomorphism $G \recht \cG$.
\end{theorem}

The cocycle superrigidity theorem yields orbital superrigidity of generalized Bernoulli actions, see \cite{P0}. In the next section, we show that the
cocycle superrigidity theorem can be used as well to prove orbital superrigidity for quotients of generalized Bernoulli actions by compact groups. In
certain cases, we completely classify these families of quotients of generalized Bernoulli actions. Note that these quotients of Bernoulli actions do
not satisfy the cocycle superrigidity theorem, but their $1$-cocycles can be completely described (see Lemma \ref{lem.one-cohom}).

\section{Superrigidity for quotients of generalized Bernoulli actions} \label{sec.orbit-equivalence-superrigid}

Theorem \ref{thm.orbital-superrigid} stated above, is a corollary to the following more precise statement.

\begin{theorem} \label{thm.superrigid}
Let $G$ be a countable group and $K$ a compact group. Let $G \times K$
act on $(X,\mu)$. We suppose that $G \actson (X,\mu)$ is essentially
free, weakly mixing and satisfies the
conclusions of the cocycle superrigidity theorem. Assume that $G \actson
X^K$ remains essentially free.

If $\Gamma \actson Y$ is any essentially free ergodic action and $\pi : X^K \recht Y$ a stable orbit equivalence, then there exists
\begin{enumerate}
\item a closed normal subgroup $K_1 \subset K$ with finite quotient
  $K_0 := K/K_1$;
\item a homomorphism $\theta : G \times K_0 \recht \Gamma$ with image
  of finite index in $\Gamma$ and with finite kernel trivially
  intersecting $K_0$;
\item a conjugation of $\Gamma$-actions
$$\Delta : Y \recht \Ind_\theta^\Gamma((G \times K_0) \actson X^{K_1})
\; ;$$
\end{enumerate}
such that $\Delta \circ \pi$ is the canonical stable orbit
equivalence. In particular, the compression constant of $\pi$ is given by
$$c(\pi) = \frac{|K_0| [\Gamma : \Im \theta]}{|\Ker \theta|} \; .$$
\end{theorem}
\begin{proof}
One can give a proof along the lines of 5.11 in \cite{P0}, using a
purely von Neumann algebraic framework. Instead, we use an approach
similar to 4.7 in \cite{V}, using the measure theoretic framework
developed in \cite{Fur2}, Theorem 3.3.

Extend $\pi$ to an orbit preserving measurable map $p : X^K \recht Y$
with associated Zimmer $1$-cocycle $\alpha : G \times X^K \recht \Gamma$ defined by
$p(g \cdot x) = \al(g,x) \cdot p(x)$ almost everywhere.
Consider the
(infinite) measure space $X^K \times \Gamma$ with the measure preserving action of $G \times \Gamma$ given by
$$g \cdot (x,s) \cdot t = (g \cdot x, \al(g,x)st) \; .$$
It follows from Theorem 3.3 in \cite{Fur2} that the $G$-action on $X^K
\times \Gamma$ admits a fundamental domain of finite measure and that the actions
$\Gamma \actson Y$ and $\Gamma \actson G \backslash (X^K \times \Gamma)$ are conjugate.

From now on, write $x \mapsto \ox$ for the quotient map $X \recht
X^K$. Denote the action of $k \in K$ on $x \in X$ as $x \cdot k$. By
our assumption, there exists a measurable map $w : X \recht \Gamma$ and a homomorphism $\theta_1 : G \recht \Gamma$ such that
$$\al(g,\ox) = w(g \cdot x) \theta_1(g) w(x)^{-1}$$
almost everywhere. Fix $k \in K$. Writing $F : X \recht \Gamma : F(x)
= w(x)^{-1} w(x \cdot k)$, we get
$$F(g \cdot x) = \theta_1(g) F(x) \theta_1(g)^{-1}$$
almost everywhere. Weak mixing implies that $F$ takes constantly the
value $s \in \Gamma$ with $s$ centralizing $\theta_1(G)$. Denote $s=
\theta_2(k)$. We have
found a continuous homomorphism $\theta_2 : K \recht \Gamma$
satisfying
\begin{equation}\label{eq.hulp}
w(x \cdot k) = w(x) \theta_2(k)
\end{equation}
almost
everywhere. Moreover, $\theta_1(G)$ and $\theta_2(K)$ commute.

Set $K_1 = \Ker \theta_2$. By compactness, $K_0 := K/K_1$ is finite and we write $\theta : G \times K_0 \recht \Gamma : \theta(g,\ok) = \theta_1(g)
\theta_2(k)$. Observe that the essential freeness of $G \actson X^K$
together with~\eqref{eq.hulp}, yields the essential freeness of $(G \times K_0) \actson X^{K_1}$.
Consider the commuting actions of $G \times K_0$ and $\Gamma$ on $X^{K_1} \times \Gamma$ given by
\begin{equation}\label{eq.hulptwee}
(g,\ok) \cdot (\ox,s) \cdot t = (\overline{g \cdot x \cdot
  k^{-1}},\al(g,\ox) s t) \; .
\end{equation}
Let $\psi : X^{K_1} \times \Gamma \recht X^K \times \Gamma$ be the quotient
  map in the first variable. It follows that the action of $G \times
  K_0$ on $X^{K_1} \times \Gamma$ admits a fundamental domain of
  finite measure and that
 $\psi$ yields a conjugacy of the actions $\Gamma
  \actson G \backslash (X^K \times \Gamma)$ and $\Gamma \actson (G \times K_0)
  \backslash (X^{K_1} \times \Gamma)$.

Finally, we have $w : X^{K_1} \recht \Gamma$ and the map $(\ox,s) \mapsto (\ox, w(x)^{-1}s)$ conjugates the actions in~\eqref{eq.hulptwee} with the
actions given by
$$(g,\ok) \cdot (\ox,s) \cdot t = (\overline{g \cdot x \cdot
  k^{-1}},\theta(g,\ok) s t) \; .$$
Since we know that this action of $G \times K_0$ admits a fundamental domain of
  finite measure, $\Ker \theta$ must be finite and $\Im \theta$ 
  of finite index. Moreover, the action $\Gamma \actson (G \times K_0)
  \backslash (X^{K_1} \times \Gamma)$ is precisely $\Gamma \actson \Ind_\theta^\Gamma((G
  \times K_0) \actson X^{K_1})$. The preceding paragraphs then yield
  that the latter is conjugate to the action $\Gamma \actson Y$. Hence, we are done.
\end{proof}

\begin{corollary} \label{cor.orbital-superrigid}
Under the conditions of the previous theorem, if $G$ does not have
finite normal subgroups, any essentially free ergodic action that is orbit equivalent with $G
\actson X^K$ is conjugate to the latter.
\end{corollary}

\section{Classification results and computations of Out} \label{sec.orbit-equivalence-classif}

By Theorem \ref{thm.superrigid}, stable orbit equivalences between quotients of generalized Bernoulli actions of weakly rigid groups are reduced to
conjugations of these actions. In this section, we prove that under condition $\cB$ (Def.\ \ref{def.conditionB}), such conjugation necessarily has a
very specific form. In that way, we prove Theorem \ref{thm.classification} stated above.

\begin{definition} \label{def.delta-conjugation}
Let $G \actson (X,\mu)$ and $\Gamma \actson (Y,\eta)$. A \emph{$\delta$-conjugation} of these actions consists of a measure space isomorphism $\Delta
: X \recht Y$ and a group isomorphism $\delta : G \recht \Gamma$ satisfying
$$\Delta(g \cdot x) = \delta(g) \cdot \Delta(x) \quad\text{almost
  everywhere.}$$
We denote by $\Aut^*(G \actson X)$ the group of all measure space
  isomorphisms $\Delta : X \recht X$ for which there exists a $\delta$
  such that $\Delta$ is a $\delta$-conjugation.

We use as well the terminology of $\delta$-conjugations for groups acting
on sets.
\end{definition}

The proof of Theorem \ref{thm.classification} consists of two steps. In Lemma \ref{lemma.lift}, we show that a conjugation of quotient actions comes
essentially from a conjugation of the original actions. Next, it is shown that, under condition $\cB$ (Def.\ \ref{def.conditionB}), a conjugation of
generalized Bernoulli actions comes from an isomorphism between the base spaces and a reshuffling of the index set.

\begin{lemma} \label{lemma.lift}
Let $G,\Gamma$ be countable and $K,L$ compact. Suppose that $(G \times
K) \actson (X,\mu)$ and $(\Gamma \times L) \actson (Y,\eta)$. Suppose
moreover that
\begin{itemize}
\item the actions $G \actson (X,\mu)$ and $\Gamma \actson (Y,\eta)$
  are weakly mixing and satisfy the conclusion of the cocycle superrigidity theorem
  \ref{thm.cocycle-superrigid};
\item there exists a $K$-equivariant measurable map $u : X \recht K$
  and an $L$-equivariant measurable map $Y \recht L$.
\end{itemize}
If $\delta : G \recht \Gamma$ is a group isomorphism and $\oDelta :
X^K \recht Y^L$ a $\delta$-conjugation of the actions $G \actson X^K$
and $\Gamma \actson Y^L$, there exists a measure space isomorphism
$\Delta : X \recht Y$, a group isomorphism $\theta : K \recht L$ and a
homomorphism $\al : G \recht \cZ(L)$ of $G$ to the center of $L$
such that
\begin{itemize}
\item $\Delta$ is a $\rho$-conjugation of
  the actions $(G \times K) \actson X$ and $(\Gamma \times L) \actson
  Y$, where $\rho$ is the isomorphism defined by $\rho(g,k) = (\delta(g),\al(g)\theta(k))$;
\item $\oDelta(\ox) = \overline{\Delta(x)}$ for almost all $x \in X$.
\end{itemize}
\end{lemma}

\begin{proof}
Let $u : X \recht K$ be a $K$-equivariant map. We define the \lq universal\rq\
$1$-cocycle $\om_G$ for the action $G \actson X^K$ with values in $G \times K$ given by
$$\om_G(g,\ox) = (g,u(g \cdot x)u(x)^{-1}) \; .$$
By Lemma \ref{lem.one-cohom} following this proof,
every $1$-cocycle for $G \actson X^K$ with values in a finite type Polish group $\cG$, is cohomologous with $\theta \circ \om_G$, for some $\theta
\in \Hom(G \times K,\cG)$, uniquely determined up to conjugacy by an element of $\cG$.

We analogously take $v : Y \recht L$ and the universal $1$-cocycle
$\om_\Gamma$. Take $\theta_1 : G \times K \recht \Gamma \times L$ and
$\theta_2 : \Gamma \times L \recht G \times K$ such that
\begin{align*}
\om_\Gamma \circ (\delta \times \oDelta) & \sim \theta_1 \circ \om_G
\; , \\
\om_G & \sim \theta_2 \circ \om_\Gamma \circ (\delta \times \oDelta)
\end{align*}
It follows that $\theta_1 \circ \theta_2$ and $\theta_2 \circ
\theta_1$ are both inner automorphisms. In particular, $\theta_1$
is a group isomorphism. If $p : \Gamma \times L \recht \Gamma$ is the
projection homomorphism, it also follows that $(p \circ \theta_1)
\circ \om_G$ is cohomologous with $(g,\ox) \mapsto \delta(g)$. Hence, $p
\circ \theta_1 = (\Ad s) \circ \delta$ for some $s \in
\Gamma$. Changing $\theta_1$, we may assume that $p \circ \theta_1 =
\delta$. In particular, $\theta_1(g,e) = (\delta(g),\al(g))$, where
$\al : G \recht L$ is a homomorphism. Also, $\theta_1(\{e\} \times K)
\subset \{e\} \times L$. Using $\theta_2$, this inclusion is an
equality, yielding an isomorphism $\theta : K \recht L$ such that
$\theta_1(g,k) = (\delta(g),\al(g) \theta(k))$. It follows that $\al$
takes values in the center $\cZ(L)$ of $L$.

Changing the equivariant map $v : Y \recht L$, we can assume
that $\om_\Gamma \circ (\delta \times \oDelta) = \theta_1 \circ \om_G$. Writing the measure space isomorphisms
$$\vphi : X \recht X^K \times K : \vphi(x) = (\ox,u(x))
\quad\text{and}\quad \psi : Y \recht Y^L \times L : \psi(y) = (\overline{y},v(y)) \; ,$$ we set $\Delta = \psi^{-1} \circ (\oDelta \times \theta)
\circ \vphi$, yielding a measure space isomorphism $\Delta : X \recht Y$. Since $\psi^{-1}(\overline{y},l) = y \cdot (v(y)^{-1} l)$, one checks
easily that $\Delta$ is a $\rho$-conjugation of $(G \times K) \actson
X$ and $(\Gamma \times L) \actson Y$, where $\rho(g,k) = (\delta(g),\al(g)\theta(k))$.
\end{proof}

We made use of the following lemma.

\begin{lemma} \label{lem.one-cohom}
Let $G$ be a countable group and $K$ a compact group. Let $G \times K$
act on $(X,\mu)$ and suppose that the action $G \actson (X,\mu)$ is
weakly mixing and
satisfies the conclusion of the cocycle superrigidity theorem \ref{thm.cocycle-superrigid}.
\begin{enumerate}
\item Denote by $\Hom_\si(G \times K,\cG)$ the set of continuous homomorphisms
$\theta : G \times K \recht \cG$ for which there exists a measurable
map $u : X \recht \cG$ such that $u(x \cdot k)
  = u(x) \theta(k)$ for almost all $x \in X$, $k \in K$. Define, for
  $\theta \in \Hom_\si(G \times K,\cG)$,
$$\om_\theta \in Z^1(G \actson X^K, \cG) : \om_\theta(g,\ox) = u(g
\cdot x) \theta(g) u(x)^{-1}$$
and note that the cohomology class of $\om_\theta$ does not depend on
the choice of $u$.
\item If $\om \in Z^1(G \actson X^K, \cG)$, there exists $\theta \in
  \Hom_\si(G \times K,\cG)$, unique up to conjugacy by an element of
  $\cG$, such that $\om \sim \om_\theta$.
\end{enumerate}
\end{lemma}
\begin{proof}
It is straightforward to check that $\om_\theta \in Z^1(G \actson X^K,
\cG)$. Using weak mixing, one gets that $\om_{\theta_1} \sim \om_{\theta_2}$ if and only if there
exists $s \in \cG$ such that $\theta_2 = (\Ad s) \circ \theta_1$.

So, take $\om \in Z^1(G \actson X^K, \cG)$. By assumption, we can take
$u : X \recht \cG$ and $\theta_1 : G \recht \cG$ such that
$$\om(g,\ox) = u(g \cdot x) \theta_1(g) u(x)^{-1} \; .$$
Fix $k \in K$ and set $F : X \recht \cG : F(x) = u(x)^{-1} u(x \cdot
k)$. Then, $F(g \cdot x) = \theta_1(g) F(x) \theta_1(g)^{-1}$. Weak
mixing implies that $F$ takes a constant value commuting with
$\theta_1(G)$ and we denote it by $\theta_2(k)$. Writing $\theta(g,k)
= \theta_1(g) \theta_2(k)$, we have $\om = \om_\theta$.
\end{proof}

Theorem \ref{thm.classification} stated in the introduction is a consequence of Theorem
\ref{thm.orbital-superrigid} and the following statement (see Remark
\ref{rem.final} below).

\begin{theorem} \label{thm.classification-version}
For every inclusion of groups $G_0 \subset G$ and every faithful
action $K \actson (X_0,\mu_0)$ of a compact group $K$, we
denote by $\si(G_0 \subset G, K \actson (X_0,\mu_0))$ the quotient of the generalized Bernoulli action given by
$$G \actson \Bigl(\prod_{G/G_0} (X_0,\mu_0)\Bigr)^K \; .$$
Denote by $\cR(G_0 \subset G,K \actson (X_0,\mu_0))$ the equivalence
relation given by the $G$-orbits.

Let $G_0 \subset G$ and $\Gamma_0 \subset \Gamma$ satisfy condition $\cB$ of Definition \ref{def.conditionB}. Let $\pi$ be a stable orbit equivalence
between $\si(G_0 \subset G, K \actson (X_0,\mu_0))$ and $\si(\Gamma_0
\subset \Gamma, L \actson (Y_0,\eta_0))$. Then, the compression
constant of $\pi$ equals $1$ and there exists
\begin{itemize}
\item an isomorphism $\delta : G \recht \Gamma$ satisfying
  $\delta(G_0) = \Gamma_0$;
\item a measure space isomorphism $\Delta_0 : (X_0,\mu_0) \recht
  (Y_0,\eta_0)$ conjugating the actions of $K$ and $L$ through an
  isomorphism $\theta : K \recht L$;
\item a homomorphism $\al : G \recht \cZ(L)$ of $G$ to the center of
  $L$ satisfying $\al(G_0) = \{e\}$;
\item an element $\psi$ of the full group of the equivalence relation
  $\cR(\Gamma_0 \subset \Gamma, L \actson (Y_0,\eta_0))$;
\end{itemize}
such that
$$(\psi \circ \pi)(\ox) = \overline{\Delta(x)} \quad\text{where}\quad
\Delta : X \recht Y : \Delta(x)_{\overline{\delta(g)}} = \al(g) \cdot \Delta_0(x_{\overline{g}})
\; .$$
\end{theorem}

\begin{proof}
We first apply Theorem \ref{thm.superrigid}. By absence of finite normal
subgroups and by weak mixing of the quotients of generalized
Bernoulli actions, we get a group isomorphism $\delta : G \recht \Gamma$ and
an element $\psi$ in the full group of the equivalence relation
  $\cR(\Gamma_0 \subset \Gamma, L \actson (Y_0,\eta_0))$ such that
  $\oDelta := \psi \circ \pi$ is a $\delta$-conjugation of the actions $G \actson
  X^K$ and $\Gamma \actson Y^L$.

The cocycle superrigidity theorem \ref{thm.cocycle-superrigid} and
Lemma \ref{lem.freeness} allow to apply Lemma
\ref{lemma.lift}. We get a measure space isomorphism $\Delta : X \recht
Y$, a group isomorphism $\theta : K \recht L$ and a homomorphism $\al
: G \recht \cZ(L)$ such that $\Delta$ is
a $\rho$-conjugation and such that $\oDelta(\ox) =
\overline{\Delta(x)}$ almost everywhere. Here $\rho(g,k) =
(\delta(g),\al(g)\theta(k))$.

Write $A_0 = L^\infty(X_0,\mu_0)$ and $B_0 = L^\infty(Y_0,\eta_0)$. Set $A = \bigotimes_{x \in G/G_0} A_0$ and $B = \bigotimes_{y \in
\Gamma/\Gamma_0} B_0$. For $x \in G/G_0$ and $y \in \Gamma/\Gamma_0$, we have natural embeddings $\pi_x : A_0 \recht A$ and $\pi_y : B_0 \recht B$.
Denote $\theta : A \recht B$ the $\rho$-conjugation given by $\theta(F) = F \circ \Delta^{-1}$.

Writing the homomorphism $\be : G \recht \cZ(K) : \be(g) = \theta^{-1}(\alpha(g)^{-1})$, observe that $\theta$ conjugates the action
$(g,\beta(g))_{g \in G}$ with the action $(\delta(g),e)_{g \in G}$. The subalgebra $\pi_{\overline{e}}(A_0)$ is quasi-periodic under the action $(g,\beta(g))_{g \in G_0}$. Hence, $\theta(\pi_{\overline{e}}(A_0))$ is quasi-periodic under the action of $\delta(G_0) \times \{e\}$. Hence, the latter is not
weakly mixing and we find $y \in \Gamma/\Gamma_0$ with $\delta(G_0) \cdot y$ finite. Composing $\Delta$ with the automorphism given by an element of
$\Gamma$, we may assume that $y={\overline{e}}$ in $\Gamma/\Gamma_0$. Observe that $\Stab {\overline{e}} = \Gamma_0$.

Write $\delta(G_0) \cap \Gamma_0 = \delta(G_1)$, where $G_1$ is a finite index subgroup of $G_0$. Whenever $x \in G/G_0$ and $x \neq {\overline{e}}$,
the orbit $G_0 \cdot x$ is infinite and hence, the orbit $G_1 \cdot x$
as well. Weak mixing implies that the fixed points of $(g,\beta(g))_{g
  \in G_1}$ are exactly $\pi_{\overline{e}}(A_0)^{\beta(G_1)}$. So,
$$\pi_{\overline{e}}(B_0) \subset B^{\delta(G_1) \times \{e\}} =
\theta(\pi_{\overline{e}}(A_0)^{\beta(G_1)}) \; .$$ It follows that $\pi_{\overline{e}}(A_0)^{\beta(G_1)}$ has a subspace that has dimension at least
$2$ and that is pointwise invariant under $\rho^{-1}(\Gamma_0 \times \{e\})$. Since for all $(g,k)$ with $g \not\in G_0$ and all $a \in
\pi_{\overline{e}}(A_0)$ with $\tau(a) = 0$, we have $\si_{(g,k)}(a)$
orthogonal to $a$, we conclude that $\rho^{-1}(\Gamma_0 \times \{e\})
\subset G_0 \times K$. Hence,
$\delta^{-1}(\Gamma_0) \subset G_0$. We
analogously get the converse inclusion. So, $\delta(G_0) = \Gamma_0$. It follows that $\pi_{\overline{e}}(B_0) \subset
\theta(\pi_{\overline{e}}(A_0)^{\beta(G_0)})$ and hence
$$\theta(A) = B \subset \theta\bigl(\underset{G/G_0}{\otimes}
\pi_{\overline{e}}(A_0)^{\beta(G_0)} \bigr) \; .$$ We conclude that $\beta(G_0) = \{e\}$. So, $\theta : \pi_{\overline{e}}(A_0) \recht
\pi_{\overline{e}}(B_0)$ is an isomorphism, implemented by a measure space isomorphism $\Delta_0 : (X_0,\mu_0) \recht (Y_0,\eta_0)$.

It follows that $\Delta_0$ is a $\theta$-conjugation of the actions $K \actson (X_0,\mu_0)$ and $L \actson (Y_0,\eta_0)$. By construction, $\Delta$
is given by the formula $\Delta(x)_{\overline{\delta(g)}} = \al(g) \cdot \Delta_0(x_{\overline{g}})$.
\end{proof}

\begin{remark} \label{rem.final}
Let $\cR:=\cR(G_0 \subset G, K \actson (X_0,\mu_0))$.
In order to deduce Theorem \ref{thm.classification} from Theorem
\ref{thm.classification-version} above, it remains to check that $\Inn
\cR$ is closed in $\Aut \cR$ and that the formula for $\Out \cR$
holds. We have the homomorphisms
$$\Hom(G/G_0 \recht \cZ(K)) \rtimes \Bigl( \frac{\Aut(G_0 \subset G)}{\Ad G_0} \times
 \frac{\Aut^*(K \actson
  X_0)}{K}\Bigr) \;\;\overset{\displaystyle\eps}{\longrightarrow}\;\;
 \frac{\Aut^*(G \actson X^K)}{G} \;\;\recht\;\; \Out(\cR) \; .$$
By Theorem \ref{thm.classification-version} the composition of both homomorphisms is
surjective. By weak mixing, the second homomorphism is injective. So, it remains to show that $\eps$ is injective. Composing with the natural
homomorphism $\frac{\Aut^*(G \actson X^K)}{G} \recht \Out(G)$, it follows that the kernel of $\eps$ is included in $\Hom(G/G_0 \recht \cZ(K)) \rtimes
\frac{\Aut^*(K \actson X_0)}{K}$. We claim that
$$\Hom(G/G_0 \recht \cZ(K)) \rtimes \frac{\Aut^*(K \actson X_0)}{K} \recht \Aut(X^K,\mu^K)$$
is injective, which suffices to get the formula for $\Out \cR$. So, let $\theta \in \Aut(K)$, $\Delta_0 \in \Aut^\theta(K \actson X_0)$ and $\al \in \Hom(G/G_0 \recht
\cZ(K))$. Define, for every $x \in G/G_0$, the automorphism $\Delta_x
\in \Aut^\theta(K \actson X_0)$ given by $\Delta_x(y) = \al(x) \cdot \Delta_0(y)$. Define $\Delta := \prod_{G/G_0} \Delta_x \in \Aut^\theta(K \actson X)$ and denote by $\oDelta$ its passage to the quotient $X^K$. Suppose
that $\oDelta = \id$ almost everywhere. We have to show that $\Delta_0 \in K$ and $\al(g) = e$ for all $g$.

Write $G/G_0 = I_1 \sqcup I_2$ with both $I_i$ infinite. Define $X_i := \prod_{I_i} (X_0,\mu_0)$ and $\Delta_i := \prod_{x \in I_i} \Delta_x$.
Observe that $\Delta_i \in \Aut^\theta(K \actson X_i)$. It suffices to prove the existence of $k \in K$ such that both $\Delta_1$ and $\Delta_2$ are
given by the action of the element $k$. Our assumption yields a measurable function $\vphi : X = X_1 \times X_2 \recht K$ such that $\Delta_i(x_i) =
x_i \cdot \vphi(x_1,x_2)$ for almost all $(x_1,x_2) \in X_1 \times X_2$. By Lemma \ref{lem.freeness}, the actions $K \actson X_i$ are essentially
free. This implies that $\vphi(x_1,x_2)$ only depends on $x_1$ and only depends on $x_2$, i.e.\ $\vphi$ is essentially constantly equal to $k \in K$.
This is what we wanted to prove.

Although one can prove by hand that $\Inn \cR$ is
closed in $\Aut \cR$, this follows as well from the formula for $\Out
\cR$. Indeed, writing $\cG = \Inn \cR \rtimes \Aut^*(G \actson X^K)$, we have a surjective continuous homomorphism between
Polish groups $\eta : \cG \recht \Aut \cR$. Note that the kernel of
$\eta$ equals the image of the embedding $G \recht \cG$. It follows
that $\Aut \cR$ is homeomorphic with $\cG / \Ker \eta$. Hence, $\Inn
\cR$ is closed in $\Aut \cR$ because $\Inn \cR \rtimes G$ is closed in $\cG$.
\end{remark}

In fact, combining Theorem \ref{thm.superrigid} with Lemmas
\ref{lem.freeness} and \ref{lemma.lift}, we immediately get the
following proposition.

\begin{proposition}
The family of actions
$$
G \actson \Bigl(\prod_G (K,\text{Haar}) \Bigr)^K \; ,
$$
where $G$ runs through the $w$-rigid groups without finite normal
subgroups and where $K$ runs through the non-trivial compact second
countable groups, consists of non-stably orbit equivalent
actions.
\end{proposition}

\section{Von Neumann strong rigidity for generalized Bernoulli
actions} \label{sec.vNstrong}

We prove a von Neumann algebraic version of Theorem
\ref{thm.classification}. We use the following terminology.

\begin{definition} \label{def.notembed}
Let $A,N$ be finite von Neumann algebras. We say that \emph{$A$ does not embed in $N$} if every, possibly non-unital, $^*$-homomorphism $A \recht
\M_n(\C) \ot N$ is identically zero.
\end{definition}

\begin{theorem} \label{thm.vNstrong}
Let $\Gamma_0 \subset \Gamma$ and $G_0
\subset G$. Consider the $(\Gamma \actson \Gamma/\Gamma_0)$-Bernoulli
action with base $(Y_0,\eta_0)$ on $(Y,\eta)$ and the $(G \actson G/G_0)$-Bernoulli
action with base $(X_0,\mu_0)$ on $(X,\mu)$. Suppose that $t > 0$ and
that
\begin{equation}\label{eq.iso}
\pi : L^\infty(Y,\eta) \rtimes \Gamma \recht
\bigl(L^\infty(X,\mu) \rtimes G \bigr)^t
\end{equation}
is a $^*$-isomorphism.
We make the following assumptions.
\begin{itemize}
\item $\Gamma$ is an ICC group with infinite subgroup $\Gamma_0$ such that $\Gamma_0
  \cap g\Gamma_0 g^{-1}$ has infinite index in $\Gamma_0$ whenever $g \in \Gamma,g
  \not\in \Gamma_0$.
\item $\Gamma$ has a normal subgroup $\Lambda$ with the relative
  property~(T).
\item $\cL(\Gamma)$ does not embed in $L^\infty(X,\mu) \rtimes G_0$.
\item $\cL(\Lambda)$ does not embed in $\cL(G_0)$.
\item $\cL(\Gamma_0)$ does not embed in $\cL(G_0 \cap gG_0g^{-1})$
  when $g \in G,g \not\in G_0$.
\item There exist $g_1,\ldots,g_n \in G$ such that $\bigcap_{i=1}^n
  g_i G_0 g_i^{-1}$ is finite.
\item The same conditions are satisfied when interchanging the roles
  of $(\Gamma,\Gamma_0)$ and $(G,G_0)$.
\end{itemize}
Then, $t=1$ and there exists
\begin{itemize}
\item a unitary $u \in L^\infty(X,\mu) \rtimes G$,
\item an isomorphism $\delta : \Gamma \recht G$ satisfying
  $\delta(\Gamma_0) = G_0$,
\item a character $\om \in \Char(\Gamma)$,
\item a measure space isomorphism $\Delta_0 : (X_0,\mu_0) \recht (Y_0,\eta_0)$,
\end{itemize}
such that
\begin{alignat*}{2}
(\Ad u \circ \pi)(\nu_s) &= \om(s) u_{\delta(s)} \quad & & \text{for all}\;\; s
\in \Gamma \; , \\
(\Ad u \circ \pi)(a) & = \alpha(a) \quad & &\text{for all}\;\; a \in
L^\infty(Y,\eta)\; , \;\; \text{where $\al : L^\infty(Y,\eta) \recht
  L^\infty(X,\mu)$ is given by} \\  & & &  \al(a) = a \circ \Delta
\;\;\text{and}\;\; (\Delta(x))_s = \Delta_0(x_{\delta(s)}) \;\;\text{for $x \in
  X$}.
\end{alignat*}
In particular, $M:=L^\infty(X,\mu)
\rtimes G$ has trivial fundamental group, $\Inn M$ closed in $\Aut
M$ and satisfies
$$\Out M \cong \Bigl( \Char(G)
\rtimes \frac{\Aut(G_0 \subset G)}{\Ad G_0} \Bigr) \times
\Aut(X_0,\mu_0) \; .$$
\end{theorem}

The same result holds for generalized Bernoulli actions with a non-commutative base space.

\begin{theorem} \label{thm.vNstrong-factor}
The previous theorem holds as well replacing the commutative base algebra $L^\infty(X_0,\mu_0)$ by either the hyperfinite II$_1$ factor or the matrix
algebra $\M_n(\C)$ with their normalized traces.
\end{theorem}

Before proving Theorems \ref{thm.vNstrong} and \ref{thm.vNstrong-factor}, we present three natural families of groups $G_0 \subset G$ that satisfy
all the conditions of these theorems. These families of groups are used in \ref{subsec.outer-vN} to give concrete computations of outer automorphism
groups of certain II$_1$ factors.

\begin{proposition} \label{prop.families}
If there exists $i \in \{1,2,3\}$ such that $(G,G_0)$ and
$(\Gamma,\Gamma_0)$ belong to the family $\cF_i$ introduced below,
then all the conditions of Theorem \ref{thm.vNstrong} are fulfilled.

\begin{itemize}
\item $\cF_1$ consists of the groups $G_0 \subset G$ satisfying
  condition $\cD$ (Def.\ \ref{def.conditionD}).
\item $\cF_2$ consists of the groups $G_0 \subset G$ given as
  follows. Take torsion free word hyperbolic property~(T) groups $K_0$ and
  $K$ and suppose $K_0 \subset K$. Set $G = K \times K_0$ and $G_0$
  the diagonal subgroup $G_0 = \{(s,s) \mid s \in K_0\}$.
\item $\cF_3$ consists of the groups $G_0 \subset G$ given as
  follows.
\begin{itemize}
\item Take $K_0 \subset K$ where $K$ is an ICC group with the
  property~(T) and where $K_0$ is an infinite, amenable, almost
  malnormal subgroup. Suppose that $K$ contains two infinite commuting
  subgroups, one of them being non-amenable.
\item Take $L_0 \subset L$ where $L$ is an ICC group and where $L_0$ is a
  non-amenable group in the class $\cC$ of Ozawa
  \cite{OzKurosh}. Suppose that $gL_0 g^{-1} \cap L_0$ is amenable
  whenever $g \in L,g \not\in L_0$. Suppose that there exist
  $g_1,\ldots,g_n \in L$ such that $\bigcap_{i=1}^n g_i L_0 g_i^{-1}$
  is finite.
\end{itemize}
Set $G = K \times L$ and $G_0 = K_0 \times L_0$.
\end{itemize}
\end{proposition}

Note that the family $\cF_3$ in is quite large. The conditions on $K_0 \subset K$ and $L_0 \subset L$ are independent of each other. The conditions
on $L_0 \subset L$ can be realized easily, see e.g.\ Proposition \ref{prop.everygroup} below.

\begin{proof}
\mbox{}\\
{\bf Family $\cF_1$.} The only non-trivial point is to observe that
$\cL(\Lambda)$ does not embed in $L^\infty(X,\mu) \rtimes G_0$ when
$G_0$ has the Haagerup property and $\Lambda$ is an infinite property
(T) group.

{\bf Family $\cF_2$.} Let $G= K \times K_0$ with diagonal subgroup
$G_0 \cong K_0$ and let $\Gamma = L \times L_0$ with diagonal subgroup
$\Gamma_0 \cong L_0$, all satisfying the above conditions.
Note that the centralizer $C_K(k)$ of an element
in $k \in K$ is cyclic (\cite{GdlH}, Th{\'e}or{\`e}me 8.34). In particular, $G$ is
an ICC group. If $g \in G$, $g \not\in G_0$, the intersection $G_0
\cap gG_0 g^{-1}$ is given by the elements of $K_0$ centralizing a
non-trivial element of $K$. This means that $G_0
\cap gG_0 g^{-1}$ is cyclic. In particular, $G_0
\cap gG_0 g^{-1}$ has infinite index in $G_0$ and also $\cL(\Gamma_0)$
does not embed in $\cL(G_0 \cap g G_0 g^{-1})$. Note that
$\cL(\Gamma)$ does not embed in $L^\infty(X,\mu) \rtimes G_0$ because
the latter is semi-solid by \cite{OzKurosh}. Finally, taking $g_1 =
e$, $g_2 = (k,e)$ and $g_3 = (h,e)$, where $h,k$ generate a free group
in $K$, it follows that $\bigcap_{i=1}^3 g_i G_0 g_i^{-1}$ is trivial.

{\bf Family $\cF_3$.} Let $G=K \times L$, $G_0 = K_0 \times L_0$. Let
$\Gamma_0 \subset \Gamma$ be of the same kind. Since $\cL(\Gamma)$ is
the tensor product of two non-injective factors, it follows from
\cite{OzKurosh} that $\cL(\Gamma)$ does not embed in $L^\infty(X,\mu)
\rtimes G_0 = (L^\infty(X,\mu) \rtimes K_0) \rtimes L_0$, the latter
being the crossed product of an injective finite von Neumann algebra
with a group in the class $\cC$ of Ozawa \cite{OzKurosh}.

Also, suppose that $\Lambda$ is a property~(T) group that contains two infinite
commuting subgroups, one of them being non-amenable. If $\cL(\Lambda)$
embeds in $\cL(K_0 \times L_0)$, property~(T) of $\Lambda$ combined
with amenability of $K_0$ implies that $\cL(\Lambda)$
actually embeds in $\cL(L_0)$, which is impossible because the latter
is solid \cite{OzSolid}.

Note that, for $g \in G \setmin G_0$, the intersection $G_0 \cap g G_0 g^{-1}$ is either amenable, either contains
$L_0$ as a finite index subgroup. So, it
belongs to the class $\cC$ of
Ozawa. Hence, $\cL(\Gamma_0)$ does not
embed in $\cL(G_0 \cap g G_0 g^{-1})$, by solidity of the latter.
\end{proof}

Theorem \ref{thm.vNstrong} is shown using the deformation/rigidity techniques developed by the first author in \cite{P1} and \cite{P2}. We make use
of another technique due to the first author, yielding unitary intertwining of subalgebras of a II$_1$ factor by using bimodules (see \cite{P5}).
This is reviewed briefly in Subsection \ref{subsec.intertwining} below. In particular, we use the notation $A \embed{M} B$ and $A \notembed{M} B$
introduced in Definition \ref{def.intertwining}. If $A,N \subset M$ and if $A$ does not embed in $N$ (in the sense of Definition \ref{def.notembed}),
then $A \notembed{M} N$.

\begin{theorem} \label{thm.malleability}
Let $G$ be an ICC group with subgroup $G_0$. Let $G \actson (X,\mu)$ be the $(G \actson G/G_0)$-Bernoulli action with base $(X_0,\mu_0)$. Write $M =
L^\infty(X,\mu) \rtimes G$. Let $Q \subset M$ be an inclusion with the
relative property~(T). Denote by $P$ the quasi-normalizer of $Q$ in
$M$ (see Subsection \ref{subsec.intertwining}).
Assume that
\begin{itemize}
\item $P$ does not embed in $L^\infty(X,\mu) \rtimes G_0$,
\item $Q$ does not embed in $\cL(G_0)$.
\end{itemize}
Then, there exists a unitary $u \in M$ such that $uPu^* \subset \cL(G)$.
\end{theorem}
\begin{proof}
The proof is entirely analogous to 4.1 and 4.4 in \cite{P1} (see also 6.4 in \cite{V}). The mixing property used in the cited proofs is replaced by
weak mixing applying Propositions \ref{prop.weak-mixing} and \ref{prop.ex-weak-mixing}.
\end{proof}

\begin{notations}
We fix the following notations and data.
\begin{itemize}
\item When $M$ is a von Neumann algebra and $n \in \N$, we denote by
  $M^n$ the von Neumann algebra $\M_n(\C) \ot M$.
\item We fix an ICC group $G$ with infinite subgroup $G_0
  \subset G$ satisfying $G_0 \cap gG_0 g^{-1}$ of infinite index in
  $G_0$ whenever $g \in G,g \not\in G_0$. We fix $\Gamma_0 \subset
  \Gamma$ satisfying the same properties.
\item We denote by $G \actson (X,\mu)$ the $(G \actson
  G/G_0)$-Bernoulli action with base $(X_0,\mu_0)$ and analogously
  $\Gamma \actson (Y,\eta)$.
\item We write $A=L^\infty(X,\mu)$ and $B=L^\infty(Y,\eta)$. We set
$A_0 = L^\infty(X_0,\mu_0)$ and $B_0 = L^\infty(Y_0,\eta_0)$. We have, for $x \in G/G_0$ and $y \in \Gamma/\Gamma_0$, the homomorphisms $\pi_x : A_0
\recht A$, $\pi_y : B_0 \recht B$.
\end{itemize}
\end{notations}

\begin{lemma} \label{lem.one}
Assume that $\Gamma$ has the relative property~(T) with respect to the normal subgroup $\Lambda$ such that
\begin{itemize}
\item $\cL(\Gamma)$ does not embed in $A \rtimes G_0$,
\item $\cL(\Lambda)$ does not embed in $\cL(G_0)$.
\end{itemize}
Assume that the same conditions are fulfilled when interchanging the
roles of $(G_0,G)$ and $(\Gamma_0,\Gamma)$.

If for $t > 0$, $B \rtimes \Gamma = (A \rtimes G)^t$, there exists a unitary $u \in
(A \rtimes G)^t$ such that $u \cL(\Gamma) u^* = \cL(G)^t$.
\end{lemma}
\begin{proof}
By Theorem \ref{thm.malleability}, we find unitaries $u,v \in (A \rtimes G)^t$ such that $u \cL(\Gamma) u^* \subset \cL(G)^t$ and $v \cL(G)^t v^*
\subset \cL(\Gamma)$. But then, $vu \cL(\Gamma) u^*v^* \subset v
\cL(G)^t v^* \subset \cL(\Gamma)$.
By Propositions \ref{prop.weak-mixing},
\ref{prop.ex-weak-mixing}, we get $vu \in \cL(\Gamma)$ and hence the equality $u \cL(\Gamma) u^* = \cL(G)^t$ holds.
\end{proof}

\begin{lemma} \label{lem.two}
Assume that $\cL(\Gamma_0)$ does not embed in $\cL(G_0 \cap g G_0
g^{-1})$ when $g \in G$, $g \not\in G_0$. Let $p \in \cL(G)^n$ and
suppose that
$$B \rtimes \Gamma = p (A \rtimes G)^n p \quad\text{with}\quad
\cL(\Gamma) = p \cL(G)^n p \; .$$
There exists $v \in p(\M_{n,\infty}(\C) \ot \cL(G))$ satisfying $vv^*
= p$, $q:= v^* v \in \M_\infty(\C) \ot \cL(G_0)$ and
$$v^* \cL(\Gamma_0) v = q \cL(G_0)^\infty q \; .$$
We use the notation $\M_\infty(\C) = \B(\ell^2(\N))$.
\end{lemma}
\begin{proof}
Note that the relative commutant of
$\cL(\Gamma_0)$ in $\cL(\Gamma)$ equals the center of $\cL(\Gamma_0)$.

{\bf Claim.} For any non-zero central projection $z \in \cL(\Gamma_0)$, there exists a
projection $q \in \cL(G_0)^k$, a unital $^*$-homomorphism $\rho : \cL(\Gamma_0) \recht
q \cL(G_0)^k q$ and a non-zero partial isometry $v \in
z(\M_{n,k}(\C) \ot \cL(G))q$ satisfying, $v^*v = q$ and
$$x v = v \rho(x) \quad\text{for all}\quad x \in \cL(\Gamma_0) \; .$$
Note that $vv^*$ belongs to the center of $\cL(\Gamma_0)$.

We postpone the proof of the claim till the end and finish the argument.
A maximality argument yields $v \in
p(\M_{n,\infty}(\C) \ot \cL(G))$ satisfying $vv^* = p = 1_{\cL(\Gamma)}$ and
$$v^* \cL(\Gamma_0) v \subset \cL(G_0)^\infty \; .$$
Write $q = v^*v \in \cL(G_0)^\infty$. Write more explicitly the
$^*$-isomorphism
$$\pi : B \rtimes \Gamma \recht q (A \rtimes G)^\infty q : \pi(x) =
v^* x v$$
satisfying $\pi(\cL(\Gamma)) = q \cL(G)^\infty q$ and
$\pi(\cL(\Gamma_0)) \subset q \cL(G_0)^\infty q$. It remains to prove
that this last inclusion is an equality. Write $M = q(A \rtimes
G)^\infty q$, $P = q \cL(G)^\infty q$ and $Q = q \cL(G_0)^\infty
q$. Since $\pi(\cL(\Gamma)) = P$, it suffices to prove that
$E_Q(\pi(\nu_s)) = 0$ for all $s \in \Gamma$, $s \not\in
\Gamma_0$. Indeed, it then follows that $E_Q(\pi(x)) =
\pi(E_{\cL(\Gamma_0)}(x))$ for all $x \in \cL(\Gamma)$, yielding $Q =
\pi(\cL(\Gamma_0))$.

Choose $s \in \Gamma$, $s \not\in \Gamma_0$. By Lemma \ref{lem.astuce}
below, we can take $t_i,r_i \in \Gamma_0$ such that
$$\|E_{\cL(\Gamma)}(x \nu_{t_i s r_i} y)\|_2 \recht  0 \quad\text{for
  all}\quad x,y \in (B \rtimes \Gamma) \ominus \cL(\Gamma) \; .$$
Writing $w_i := \pi(\nu_{t_i s r_i})$, it follows that $\|E_P(x w_i
  y)\|_2 \recht 0$ for all $x,y \in M \ominus P$. Fix $a \in A_0$ with
  $\tau(a) = 0$, but $a \neq 0$. Note that $\pi_{\overline{e}}(a)$ and $\cL(G_0)$
  commute and write $x := (1 \ot \pi_{\overline{e}}(a))q \in M \ominus P$. On the other
  hand, since $\si_g(\pi_{\overline{e}}(a)) = \pi_{\overline{g}}(a)$,
  it is easy to compute that, for all $i$,
$$E_P(x w_i x^*) = \tau(aa^*) E_Q(w_i) \; .$$
Since $\tau(aa^*) \neq 0$, we conclude that $\|E_Q(w_i)\|_2 \recht 0$ when $n \recht
\infty$. Since $w_i = \pi(\nu_{t_i}) \pi(\nu_s) \pi(\nu_{r_i})$, where
the exterior factors are unitaries in $Q$, it follows that
$E_Q(\pi(\nu_s)) = 0$, ending the proof.

It remains to prove the claim. Let $z$ be a non-zero central
projection in $\cL(\Gamma_0)$. First of all,
$$z \cL(\Gamma_0) \embed{\cL(G)^n} \cL(G_0)^n \; .$$
Indeed, if not, we get a sequence of unitaries $w_i \in
\cL(\Gamma_0)z$ satisfying
$$\|E_{\cL(G_0)^n}(x w_i y)\|_2 \recht 0 \quad\text{for all}\quad x,y
\in \cL(G)^n \; .$$
As in Propositions \ref{prop.weak-mixing}, \ref{prop.ex-weak-mixing},
it follows that any element of $z (A \rtimes G)^n z$ that commutes
with all $w_i$, belongs to $z \cL(G)^n z$. But, $z
\pi_{\overline{e}}(B_0)$ provides elements commuting with $w_i$ and
orthogonal to $z \cL(\Gamma) z = z \cL(G)^n z$, yielding a
contradiction.

So, we get a projection $q \in \cL(G_0)^k$, a non-zero
partial isometry $v \in z(\M_{n,k}(\C) \ot \cL(G))q$ and
a unital $^*$-homomorphism $\rho : \cL(\Gamma_0)z \recht q
\cL(G_0)^k q$ satisfying $x v = v \rho(x)$ for all $x \in
\cL(\Gamma_0)z$. We have to prove that $v^*v \in \cL(G_0)^k$. Write $Q
= \rho(\cL(\Gamma_0)z)$ and note that $v^*v \in Q' \cap q \cL(G)^k q$. By our assumption, $Q$ does not embed in
$\cL(G_0 \cap g G_0 g^{-1})$ when $g \in G,g \not\in
G_0$. Propositions \ref{prop.weak-mixing}, \ref{prop.ex-weak-mixing}
imply that $Q' \cap q \cL(G)^k q \subset q \cL(G_0)^k q$. So, we are done.
\end{proof}

\begin{lemma}\label{lem.astuce}
Let $\Gamma_0 \subset \Gamma$ and suppose that the quasi-normalizer of
$\Gamma_0$ in $\Gamma$ equals $\Gamma_0$. If $s \in \Gamma, s \not\in
\Gamma_0$ and if $A,B \subset \Gamma$ are finite subsets, there exist
$t,r \in \Gamma_0$ such that $t s r \not\in A \Gamma_0 B$.
\end{lemma}
\begin{proof}
Our assumption and Lemma \ref{lem.bete} imply the following: if $A,B
\subset \Gamma$ are finite such that at least one of both is disjoint
with $\Gamma_0$, there exists $t \in \Gamma_0$
such that $t \not\in A \Gamma_0 B$.

Let $A,B \subset \Gamma$ be finite. First choose $r \in \Gamma_0$ such
that $r \not\in s^{-1}\Gamma_0 B$. It follows that $B r^{-1} s^{-1}$
is disjoint with $\Gamma_0$. So, we can take $t \in \Gamma_0$ such that
$t \not\in A \Gamma_0 B r^{-1} s^{-1}$. It follows that $t s r \not\in A \Gamma_0 B$.
\end{proof}

\begin{lemma} \label{lem.three}
Let $q$ be a projection in $\cL(G_0)^\infty$ and assume that
$$B \rtimes \Gamma = q (A \rtimes G)^\infty q \quad\text{with}\quad
\cL(\Gamma) = q \cL(G)^\infty q \quad\text{and}\quad \cL(\Gamma_0) = q
\cL(G_0)^\infty q \; .$$
Then, $B \subset q(A \rtimes G_0)^\infty q$.
\end{lemma}
\begin{proof}
Let $H=\ell^2(\N)$ and consider $(A \rtimes G)^\infty$ as $\B(H) \ot
(A \rtimes G)$.
Since $\pi_{\overline{e}}(B_0)$ commutes with $\cL(\Gamma_0)$, we find
$$\pi_{\overline{e}}(B_0) \subset \B(H) \ot
\overline{\pi_{\overline{e}}(A_0) \; \cL(G_0)} \; .$$
We also get, for all $s \in \Gamma \setmin \Gamma_0$,
$$\nu_s  \in \B(H) \ot \overline{\lspan\{u_g \mid g \in G
  \setmin G_0 \}} \; .$$
Combining both, it follows that
$$\pi_{\overline{s}}(B_0) \subset \B(H) \ot \overline{A_{I'} \cL(G)}$$
for all $s \in \Gamma - \Gamma_0$, where we write $I'= G/G_0 \setmin
\{\overline{e}\}$ and $A_{I'} = \underset{x \in I'}{\otimes} A_0$.

Let now $a \in \pi_{\overline{s}}(B_0)$ for some $s \in \Gamma - \Gamma_0$. We shall prove that $a  \in \B(H) \ot (A
\rtimes G_0)$. Let $E : \B(H) \ot (A \rtimes
G) \recht \B(H) \ot (A \rtimes G_0)$ be the natural conditional
expectation. Set $b = a - E(a)$. We argue that $b=0$. To do so, take
$x \in \pi_{\overline{e}}(B_0)$ with $\tau(x) = 0$ and $x$
invertible.

We have $xa = ax$ and since $x \in \B(H) \ot (A \rtimes
G_0)$, the same holds when we replace $a$ by $b$. Note that $a \in
\B(H) \ot \overline{A_{I'} \; \cL(G)}$, which implies that
$$b \in \B(H) \ot \overline{A_{I'} \; \lspan\{u_g \mid g \in G \setmin
  G_0\}} \; .$$
Also, $x \in \B(H) \ot \overline{\pi_{\overline{e}}(A_0 \ominus \C)
  \; \cL(G_0)}$.
But then,
\begin{align*}
bx &\in \B(H) \ot \overline{A_{I'} \; \cL(G)} \; , \\
xb &\in \B(H) \ot \overline{\pi_{\overline{e}}(A_0 \ominus \C)
  A_{I'} \; \cL(G)} \; .
\end{align*}
It follows that $bx$ and $xb$ belong to orthogonal
subspaces. Since $bx = xb$, we get $xb = 0$. Since $x$ is
invertible, we conclude that $b=0$. We have shown that $\pi_{\overline{s}}(B_0)
\subset \B(H) \ot (A \rtimes G_0)$ for all $s \in \Gamma-\Gamma_0$. We
already knew that the same holds for $s = e$ and it follows that $B
\subset \B(H) \ot (A \rtimes G_0)$.
\end{proof}

We finally prove Theorem \ref{thm.vNstrong}.

\begin{proof}[Proof of Theorem \ref{thm.vNstrong}]
Write $M = (A \rtimes G)^t$ and identify $B \rtimes \Gamma = M$
through the isomorphism $\pi$ in~\eqref{eq.iso}. We apply Lemma
\ref{lem.one} that entitles to apply Lemma \ref{lem.two}, that in turn
allows to apply Lemma \ref{lem.three}. We end up with $v \in
\M_{1,\infty}(\C) \ot M$ satisfying $vv^*=1$, $q:=v^*v \in
\cL(G_0)^\infty$ and $v^*B v \subset q(A \rtimes G_0)^\infty q$. Since
the center of $A \rtimes G_0$ equals $A^{G_0}$ and $q \in \cL(G_0)^\infty$,
it follows that $v(1 \ot z) \neq 0$ for every non-zero central
projection $z \in A \rtimes G_0$. So, we can apply Theorem
\ref{thm.normality}. By assumption, there exist $g_1,\ldots,g_n$ such
that $\bigcap_{i=1}^n g_i G_0 g_i^{-1}$ is finite. We conclude that
$B \embed{M} A^t$. By Theorem A.1 in \cite{P5}, we have $uBu^* =
A^t$ for some unitary $u \in M$ and we may assume that $B=A^t$.

Applying Theorem \ref{thm.classification-version} (and Theorem
\ref{thm.cocycle-superrigid} to deal with the scalar $1$-cohomology),
we are done. Note that $\Inn M$ is closed in $\Aut M$ by an argument
as in Remark \ref{rem.final}.
\end{proof}

\begin{proof}[Proof of Theorem \ref{thm.vNstrong-factor}]
In the proof of Theorem \ref{thm.vNstrong}, we almost did not use the commutativity of the base algebra $L^\infty(X_0,\mu_0)$. If $(A \rtimes G)^t =
B \rtimes \Gamma$, exactly the same argument as in the proof of Theorem \ref{thm.vNstrong} yields $B \embed{M} A^t$ and, by symmetry, $A^t \embed{M}
B$. Then, Lemma 8.4 in \cite{IPP} yields a unitary $u \in M$ such that $uBu^* = A^t$. Combining Theorem 5.5 and Corollary 5.9 in \cite{P4}, we arrive
at the desired conclusion.
\end{proof}

\subsection{Cocycle crossed products}

The von Neumann strong rigidity theorem \ref{thm.vNstrong} remains valid when replacing the ordinary crossed product $L^\infty(X,\mu) \rtimes G$ by a
\emph{twisted crossed product} $L^\infty(X,\mu) \twist{\Om} G$, for some scalar $2$-cocycle $\Om$ on $G$ with coefficients in $S^1$. This general
philosophy holds true as well for all the results in
\cite{P1,P2}. This yields below a twisted version of Theorem
\ref{thm.vNstrong} and it is applied in Theorem \ref{thm.without-anti} to give examples of II$_1$ factors $M$ without anti-automorphisms, with trivial outer automorphism group and with trivial
fundamental group.

We first introduce a bit of notation related to the $2$-cohomology of
a countable group $G$. We denote by $Z^2(G,S^1)$ the abelian group of functions
$\Om : G \times G \recht S^1$ that satisfy
$$\Om(g,h) \Om(gh,k)  = \Om(g,hk) \Om(h,k) \quad\text{for all}\quad g,h,k \in G \; .$$
The elements of $Z^2(G,S^1)$ are called \emph{scalar $2$-cocycles}. Whenever $\om : G \recht S^1$ is a function, we denote by $\partial \omega$ the
$2$-cocycle given by $(\partial \omega)(g,h) = \omega(gh) \overline{\om(g) \om(h)}$. Note that $\om$ is a character if and only if $\partial \omega =
1$. The $2$-cocycles of the form $\partial \omega$ are called \emph{coboundaries} and form a subgroup of $Z^2(G,S^1)$. The quotient group is denoted
by $H^2(G,S^1)$.

Whenever $G \actson (X,\mu)$ and $\Om \in Z^2(G,S^1)$, we have a \emph{twisted crossed product} $L^\infty(X,\mu) \twist{\Om} G$, generated by a copy
of $L^\infty(X,\mu)$ and unitaries $(u_g)_{g \in G}$ satisfying $u_g F(\cdot) u_g^* = F(g^{-1} \cdot)$ and $u_g u_h = \Om(g,h) u_{gh}$.

\begin{theorem}\label{thm.vNstrong-cocycle}
Let $\Gamma_0 \subset \Gamma$ and $G_0 \subset G$. Let $\Om \in Z^2(G,S^1)$ and $\Om_\Gamma \in Z^2(\Gamma,S^1)$. Consider the same data,
constructions and conditions as in Theorem \ref{thm.vNstrong}, systematically replacing crossed products by twisted crossed products. Suppose that $t
> 0$ and that
\begin{equation*}
\pi : L^\infty(Y,\eta) \twist{\Om_\Gamma} \Gamma \recht \bigl(L^\infty(X,\mu) \twist{\Om} G \bigr)^t
\end{equation*}
is a $^*$-isomorphism. Then, $t=1$ and there exists
\begin{itemize}
\item a unitary $u \in L^\infty(X,\mu) \twist{\Om} G$,
\item an isomorphism $\delta : \Gamma \recht G$ satisfying
  $\delta(\Gamma_0) = G_0$,
\item a map $\om : \Gamma \recht S^1$ such that $\Om \circ \delta =
  \Om_\Gamma \cdot (\partial \om)$,
\item a measure space isomorphism $\Delta_0 : (X_0,\mu_0) \recht (Y_0,\eta_0)$,
\end{itemize}
such that
\begin{alignat*}{2}
(\Ad u \circ \pi)(\nu_s) &= \om(s) u_{\delta(s)} \quad & & \text{for all}\;\; s
\in \Gamma \; , \\
(\Ad u \circ \pi)(a) & = \alpha(a) \quad & &\text{for all}\;\; a \in L^\infty(Y,\eta)\; , \;\; \text{where $\al : L^\infty(Y,\eta) \recht
  L^\infty(X,\mu)$ is given by} \\  & & &  \al(a) = a \circ \Delta
\;\;\text{and}\;\; (\Delta(x))_s = \Delta_0(x_{\delta(s)}) \;\;\text{for $x \in
  X$}.
\end{alignat*}
In particular, setting $M := L^\infty(X,\mu) \twist{\Om} G$,
\begin{itemize}
\item $M$ has trivial fundamental group,
\item $M$ admits an anti-automorphism if and only if $\Om$ and $\overline{\Om}$ define the same element of
$H^2(G,S^1)$,
\item the outer automorphism group $\Out(M)$ is given by
$$\Out(M) \cong \cO \times \Aut(X_0,\mu_0) \quad\text{where}\quad e \recht \Char G \recht \cO \recht \frac{\Aut_\Om(G_0 \subset G)}{\Ad G_0} \recht
e$$ is a short exact sequence and $\Aut_\Om(G_0 \subset G)$ denotes the subgroup of $\delta \in \Aut(G)$ satisfying $\delta(G_0) = G_0$ and $\Om
\circ \delta = \Om$ in $H^2(G,S^1)$.
\end{itemize}
\end{theorem}

\subsection{Intertwining by bimodules and weak mixing techniques}\label{subsec.intertwining}

We briefly review a technique developed by the first author in order to intertwine unitarily subalgebras of a II$_1$ factor using bimodules (see
\cite{P5}). We prove a few general results that were needed above and that are, in fact, of independent interest as well.

Let $(B,\tau)$ be a von Neumann algebra with normal faithful tracial state
$\tau$. Let $H_B$ be a right Hilbert
$B$-module. We say that $H_B$ is \emph{of finite type}, if $H_B$ is
isomorphic with a sub-$B$-module of $\C^n \ot L^2(B,\tau)$.

Let $(M,\tau)$ be a von Neumann algebra with faithful tracial state
$\tau$. Let $B \subset M$ be a von Neumann subalgebra. Recall the
\emph{basic construction} $\langle M,e_B \rangle$. This is the von
Neumann algebra acting on $L^2(M,\tau)$ generated by $M$ and the
orthogonal projection onto the closure of $B$ in
$L^2(M,\tau)$. Alternatively, it consists of the operators on
$L^2(M,\tau)$ that commute with the right action of $B$. Moreover,
$\langle M,e_B \rangle$ has a canonical semi-finite trace $\Tr$
characterized by $\Tr(xe_By) = \tau(xy)$ for all $x,y \in M$.

The \emph{quasi-normalizer} of $B$ in $M$ is defined as the $^*$-algebra of
elements $x \in M$ for which there exist finite subsets
$\{y_1,\ldots,y_n\}, \{z_1,\ldots,z_m\} \subset M$ satisfying
$$B x \subset \sum_{i=1}^n y_i B \quad\text{and}\quad x B \subset
\sum_{j=1}^m B z_j \; .$$

\begin{definition}\label{def.intertwining}
Let $(M,\tau)$ be a von Neumann algebra with faithful tracial state
$\tau$. Let $A,B \subset M$ be von Neumann subalgebras. We allow that
the unit elements $1_A$, $1_B$ are non-trivial projections in $M$.
The following statements are equivalent.
\begin{enumerate}
\item $1_A L^2(M,\tau) 1_B$ admits a sub-$A$-$B$-bimodule that is of finite
  type as a $B$-module.
\item There exists $a \in A' \cap 1_A \langle M,e_B \rangle^+ 1_A$
  satisfying $0 < \Tr(a) < \infty$.
\item There is no sequence of unitaries $u_i$ in $A$ satisfying
$$\|E_B(a^* u_i b)\|_2 \recht 0 \quad\text{for all}\quad a,b \in 1_A M
1_B \; .$$
\item There exists $n \in \N$, a projection $p \in \M_n(\C) \ot B$, a
  $^*$-homomorphism $\theta : A \recht p(\M_n(\C) \ot B)p$ and a
  non-zero partial isometry $v \in 1_A (M_{1,n}(\C) \ot M) p$
  satisfying
$$x v = v \theta(x) \quad\text{for all}\quad x \in A \; .$$
\end{enumerate}
If one of the statements above is satisfied, we say that \emph{$A$
  embeds in $B$ inside $M$} and we denote this relation as
$$A \embed{M} B \; .$$
\end{definition}

In the arguments above, weak mixing plays a crucial role. We use the following terminology and results. A related notion of mixing MASA's in a II$_1$
factor has been introduced and studied in \cite{JS}. All this goes back to \cite{P6}, where mixing properties were already used to determine
normalizers of subalgebras of a II$_1$ factor.

\begin{definition}
Let $Q \subset A \subset (N,\tau)$ be an inclusion of finite von Neumann algebras. We say that \emph{$A \subset N$ is weakly mixing through $Q$} if
there exists a sequence $u_n \in \cU(Q)$ such that
$$\|E_A(x u_n y)\|_2 \recht 0 \quad\text{for all}\quad x,y \in N \ominus A \; .$$
\end{definition}

\begin{proposition} \label{prop.weak-mixing}
Let $Q \subset A \subset (N,\tau)$ such that $A \subset N$ is weakly
mixing through $Q$. Any sub-$Q$-$A$-bimodule of $L^2(N,\tau)$ which is
of finite type as an $A$-module, is contained in $L^2(A,\tau)$. In
particular, if $x \in N$ satisfies $Q x \subset \sum_{i=1}^k y_i A$
for some finite subset $\{y_1,\ldots,y_k\} \subset N$, we have $x \in A$.
In particular, $N \cap Q'
\subset A$.
\end{proposition}
\begin{proof}
The proof is entirely analogous to the one of 3.1 in \cite{P1} (see also D.4 in \cite{V}).
\end{proof}

The following proposition can be easily proved.

\begin{proposition} \label{prop.ex-weak-mixing}
We have the following examples of weakly mixing inclusions.
\begin{itemize}
\item Let $G_0 \subset G$ be countable groups and $Q \subset \cL(G_0)$. Assume that $Q \notembed{\cL(G_0)} \cL(G_0 \cap g G_0 g^{-1})$ whenever $g
\in G$, $g \not\in G_0$. Then, $\cL(G_0) \subset \cL(G)$ is weakly mixing through $Q$.
\item Let $G \actson (X,\mu)$ be the $(G \actson G/G_0)$-Bernoulli action with base $(X_0,\mu_0)$ and $P \subset \cL(G)$. If $P \notembed{\cL(G)} \cL(G_0)$, the
inclusion $\cL(G) \subset L^\infty(X,\mu) \rtimes G$ is weakly mixing through $P$.
\item Let $(Y_0,\eta_0) \recht (X_0,\mu_0)$ be a quotient map (i.e.\ $L^\infty(X_0,\mu_0) \subset
L^\infty(Y_0,\eta_0)$ is a trace preserving inclusion of von Neumann algebras). Let $G \actson (X,\mu)$, resp.\ $G \actson (Y,\eta)$, be the $(G \actson G/G_0)$-Bernoulli action with base $(X_0,\mu_0)$,
resp.\ $(Y_0,\eta_0)$. Write $M = L^\infty(X,\mu) \rtimes G$ and $\Mtil = L^\infty(Y,\eta) \rtimes G$.

If $P \subset M$ and $P \notembed{M} L^\infty(X,\mu) \rtimes G_0$, then $M \subset \Mtil$ is weakly mixing through $P$.
\end{itemize}
\end{proposition}

We present another general result. The aim is the following. Consider a crossed product $A
\rtimes G$. Let $B$ be a sufficiently regular subalgebra of $A \rtimes
G$ and suppose that $B \subset A \rtimes H$ for a subgroup $H \subset G$ that is
far from being normal. Then, $B$ embeds in $A$ inside $A \rtimes G$.

Recall that an action of a countable group $G$ on a von Neumann
algebra $A$ is said to be \emph{strictly outer} if for every $g \neq
e$ and $a \in A$ satisfying $$\al_g(b) a = a b \quad\text{for
  all}\quad b \in A \; ,$$ we
have $a=0$. Strict outerness is equivalent with $(A \rtimes G) \cap A'
= \cZ(A)$.

\begin{theorem} \label{thm.normality}
Let a countable group $G$ act strictly outerly on the finite von Neumann algebra
$(A,\tau)$. Let $H$ be a subgroup of $G$. Set $M = A \rtimes
G$. Let $B \subset M$ be a von Neumann subalgebra. Suppose that the
quasi-normalizer of $B$ in $M$ is dense in $M$.

If $B \embed{M} (A \rtimes H)z$ for every non-zero central projection
$z \in A \rtimes H$, then $B \embed{M} (A \rtimes (H \cap
gHg^{-1}))z_1$ for all $g \in G$ and all non-zero central projections
$z_1 \in A \rtimes (H \cap gHg^{-1})$.
Repeating the procedure, we get in particular that $$B \embed{M} A \rtimes \bigl(\bigcap_{i=1}^n g_i H g_i^{-1} \bigr)$$ for all
$g_1,\ldots,g_n \in G$.
\end{theorem}
\begin{proof}
Observe that the strict outerness implies that $\cZ(A \rtimes
\Gamma) = \cZ(A)^\Gamma$ for any subgroup $\Gamma \subset G$.
Denote for $n \in \N$,
\begin{equation*}
\cI_n := \Bigl\{ a \in \M_{1,n}(\C) \ot M \;\Big| \quad\parbox{11.5cm}{$a$ is a partial isometry for which there exists a, possibly non-unital,
$^*$-homo\-morphism $\al : B \recht \M_n(\C) \ot B$ satisfying $xa= a \alpha(x)$ for all $x \in B$.} \quad\Bigr\}
\end{equation*}
The coefficients of elements of $\cI_n$, $n \in \N$, linearly span a subalgebra of $M$ which is dense because the quasi-normalizer of $B$ in $M$ is
dense.

Take $g \in G$. Suppose that $$B \embed{M} (A \rtimes H)z
\quad\text{and}\quad B \notembed{M} (A \rtimes (H \cap gHg^{-1}))z_1$$
for every non-zero central projection $z \in A \rtimes H$ and for some
non-zero central projection $z_1 \in  A \rtimes (H \cap gHg^{-1})$.
Take a non-zero partial
isometry $w \in \M_{1,m}(\C) \ot M$ and a, possibly non-unital, $^*$-homomorphism $\theta : B \recht \M_m(\C) \ot (A \rtimes H)$ such that $x w = w
\theta(x)$ for all $x \in B$. We claim that we may assume that the
smallest projection $p \in M$ satisfying $w^*w \leq 1 \ot p$, is
arbitrary close to $1$. Indeed, for any projection $z \in \cZ(A)^H$, we
find $w$ with $w(1 \ot z) \neq 0$. Moreover, we can take direct sums
of $w$'s and multiply $w$ on the right with an element of $A \rtimes
H$. This proves the claim.

Take as well a sequence of unitaries $v_i \in B$ such that, with $H^g := H \cap gHg^{-1}$,
$$\|z_1 E_{A \rtimes H^g}(a^* v_i b)\|_2 \recht 0 \quad\text{for all}\quad a,b \in M \; .$$
When $Q \subset M$, we continue writing $E_Q$ for the conditional expectation $\id \ot E_Q$ of $\M_n(\C) \ot M$ onto $\M_n(\C) \ot Q$. We use the
same convention for other maps that we extend to matrix spaces.

{\bf Claim.} For all $x \in M$, we have
$$\|(1 \ot z_1) E_{A \rtimes H}((1 \ot u_g) w^* v_i x)\|_2 \recht 0 \; .$$
Whenever $L \subset G$, we denote by $p_L$ the orthogonal projection of $L^2(M)$ onto the closed linear span of $\{A u_g \mid g \in L\}$. Let $k \in
G$ be arbitrary. Choose $\eps > 0$. It suffices to show that
\begin{equation}\label{eq.weips}
\|(1 \ot z_1) E_{A \rtimes H}((1 \ot u_g) w^* v_i u_k)\|_2 < 3 \eps
\end{equation}
for $i$ sufficiently large. Take $K \subset G$ finite and $w_0 \in \M_{1,m}(\C) \ot \Im p_K$ such that $\|w - w_0\|_2 < \eps$. It follows that
\begin{alignat*}{2}
E_{A \rtimes H}((1 \ot u_g) w^* v_i u_k) &= p_H((1 \ot u_g)
\theta(v_i) w^* u_k) \approx p_H((1 \ot u_g) \theta(v_i)w_0^* u_k) & &\quad\text{with error $< \eps$} \\
&= p_{H \cap gHKk}((1 \ot u_g) \theta(v_i) w_0^* u_k) \approx p_{H
  \cap gHKk}((1 \ot u_g) w^* v_i u_k) & &\quad\text{with error $< \eps$.}
\end{alignat*}
Take $L \subset G$ finite such that $H \cap gHKk = \bigsqcup_{s
  \in L} H^g s$. Hence, we get
$$(1 \ot z_1) p_{H \cap gHKk}((1 \ot u_g) w^* v_i u_k) = \sum_{s \in
  L} (1 \ot z_1) E_{A \rtimes
  H^g}((1 \ot u_g) w^* v_i u_{ks^{-1}})u_s \; .$$
Our choice of $v_i$ implies that every term on the right hand side
  converges to $0$ in $L^2$-norm. Combined with the estimate above, we have shown~\eqref{eq.weips} for $i$ sufficiently large. This proves the claim.

Let $k \in \N$ and $a \in \cI_k$. Take $\al : B \recht \M_k(\C) \ot B$ such that $x a = a \al(x)$ for all $x \in B$. Observe that
$$v_i a (1 \ot w) = a (1 \ot w) (\id \ot \theta)\al(v_i) \; .$$
Since $(\id \ot \theta)\al(v_i) \in \M_{km}(\C) \ot (A \rtimes H)$, it follows that
$$\|(1 \ot z_1)E_{A \rtimes H}((1 \ot u_g)w^* v_i a (1 \ot w)) \|_2 =
\|(1 \ot z_1) E_{A \rtimes H}((1 \ot u_g)w^* a (1 \ot w))\|_2$$
for all $i$. Applying our claim, we conclude that
$$(1 \ot z_1) E_{A \rtimes H}((1 \ot u_g)w^* a (1 \ot w)) = 0 \quad\text{whenever}\quad a \in \cI_k \; .$$
Since the coefficients of $a, a \in \cI_k, k \geq 1$ span a strongly
dense subset of $M$, it follows that $(1 \ot
z_1) E_{A \rtimes H}((1 \ot u_g)w^* x w)=0$ for all $x \in M$. Let $q
\in M$ be the smallest projection satisfying $w^* w \leq 1 \ot q$. It
follows that $z_1 E_{A \rtimes H}(u_g qxq)=0$ for all $x \in M$. Since
$q$ can be taken arbitrarily close to $1$, we arrive at the
contradiction $z_1 = 0$.
\end{proof}

\section{Examples and computations} \label{sec.examples}

In order to illustrate Theorem \ref{thm.classification-version}, we present several examples of inclusions $G_0 \subset G$ satisfying condition $\cB$
(Def.\ \ref{def.conditionB}). This yields in Subsection \ref{subsec.outer-equivalence} continuous families of type II$_1$ equivalence relations $\cR$
with $\Out(\cR)$ an arbitrary countable group.

We present as well inclusions $G_0 \subset G$ satisfying the stronger condition $\cD$ (Def.\ \ref{def.conditionD}) and inclusions belonging to the
family $\cF_3$ of Proposition \ref{prop.families}. This yields in
Subsection \ref{subsec.outer-vN} continuous families of II$_1$ factors $M$ with $\Out(M)$ an
arbitrary finitely presented group.

\subsection{Some outer automorphism groups of discrete groups}

Let $G = \SL(n,\Z)$, $n \geq 3$. By \cite{HR} (see also Example 2.6 in \cite{Borel}), we get the following.
\begin{itemize}
\item If $n$ is odd, $\Out(G)$ has two elements, the non-trivial one
  being $\si(A) = (A^{-1})^t$.
\item If $n$ is even, $\Out(G)$ has four elements and is generated by
  the automorphism $\si$ above and the automorphism $\al(A) =
  TAT^{-1}$, where $T \in \GL(n,\Z)$ has determinant $-1$.
\item Since $G$ equals its commutator subgroup, $\Char G$ is trivial.
\end{itemize}

We give a self-contained and elementary argument for the following two, probably well known, computations of $\Out(G)$.

\begin{proposition} \label{prop.outer}
Let $n \geq 2$ and write $G = \Z^n \rtimes \GL(n,\Z)$, $G^+=\Z^n \rtimes \SL(n,\Z)$.
\begin{itemize}
\item $\Out(G)$ is trivial and $\Out(G^+)$ has two elements, a
  non-trivial one being given by conjugating by an element of $G
  \setmin G^+$.
\item If $n \geq 3$, $\Char(G) = \{1,\det\}$ and $\Char(G^+)$ is trivial.
\item $G$ and $G^+$ are ICC groups.
\end{itemize}
\end{proposition}
\begin{proof}
Let $\al \in \Aut(G)$. We claim that $\al(\Z^n) = \Z^n$. Write $\pi : G \recht \GL(n,\Z)$ the quotient map. Invoking knowledge on the normal subgroups
of $\GL(n,\Z)$, one gets $\pi(\al(\Z^n)) \subset \{\pm 1\}$ and hence $\al(\Z^n) \subset \Z^n$. The following elementary argument proves as well the
claim. Define
$$\Gamma = \{ x \in \Z^n \mid \pi(\al(x)) = 1 \} \; .$$
It is clear that $\Gamma$ is globally $\GL(n,\Z)$-invariant. Suppose $\Gamma \neq \Z^n$ and take $a \in \Z^n$ with $\al(a) = (x_0,g_0)$ and $g_0 \neq
1$. Normality yields $((1-g_0)y + x_0,g_0) \in \al(\Z^n)$ for all $y
\in \Z^n$. So, $((1-g_0)y,1) \in \al(\Z^n)$ for all $y \in \Z^n$. Since $g_0
\neq 1$, $\Gamma \neq \{0\}$. Since $\Gamma$ is globally $\GL(n,\Z)$-invariant, it follows that $\Gamma$ is of finite index in $\Z^n$. So,
$\pi(\al(\Z^n))$ is finite. Expressing that $((1-g_0)y,1)$ and
$(x_0,g_0)$ commute for all $y \in \Z^n$ (since both are elements of $\al(\Z^n)$), it follows that
$(1-g_0)^2 = 0$. Hence, $g_0 = 1 + A$ with $A \neq 0$ and $A^2=0$. But then $\pi(\al(a^k)) = 1 + k A$ implies $\pi(\al(\Z^n))$ to be infinite;
contradiction. So, $\al(\Z^n) \subset \Z^n$ and applying the same to $\al^{-1}$, we get equality.

Take $T \in \GL(n,\Z)$ such that $\al(x,1) = (Tx,1)$ for all $x \in \Z^n$. Composing with $\Ad(0,T)$, we get $\al(x,1) = (x,1)$ for all $x \in \Z^n$.
This implies that $\al(x,g) = (x+\delta(g),g)$ where $\delta : \GL(n,\Z) \recht \Z^n$ is a $1$-cocycle: $\delta(gh) = \delta(g) + g \delta(h)$. In
order to obtain $\Out(G)$ trivial, it remains to show that every such $1$-cocycle is of the form $\delta(g) = (1 - g)x$ for some $x \in \Z^n$.

Consider a $1$-cocycle $\delta : \SL(2,\Z) \recht \Z^2$. We prove that $\delta$ is a coboundary. The group $\SL(2,\Z)$ is generated by the element
$s=\bigl(\begin{smallmatrix} 0 & -1 \\ 1 & 1 \end{smallmatrix}\bigr)$ of order $6$ and the element $f=\bigl(\begin{smallmatrix} 0 & 1 \\ -1 & 0
\end{smallmatrix}\bigr)$ of order $4$. They satisfy $s^3=f^2=-1$. Since $1-s$ belongs to $\SL(2,\Z)$, we can replace $\delta$ by a cohomologous
$1$-cocycle and suppose that $\delta(s) = 0$ and hence, $\delta(s^k) = 0$ for all $k$. In particular, $\delta(f^2)=0$. Since $1+f$ is an invertible
matrix (over $\Q$), it follows that $\delta(f) = 0$. Since $s$ and $f$ generate $\SL(2,\Z)$, we get $\delta(g) = 0$ for all $g \in \SL(2,\Z)$.

Consider next a $1$-cocycle $\delta : \SL(n,\Z) \recht \Z^n$, for $n \geq 3$. Let $E_{ij}$ be the matrix with zeros everywhere and $1$ in position
$ij$. We have the commutator
\begin{equation}\label{eq.commutator}
[1+nE_{ij},1+E_{jk}] = 1+nE_{ik}
\end{equation}
whenever $i,j,k$ are distinct. Write $\delta_{ij}:= \delta(1+E_{ij})$. Applying the cocycle relation to~\eqref{eq.commutator} with $n=1$, we
get
\begin{equation}\label{eq.hulpje}
\delta_{ik} = (E_{ij} - E_{ik})\delta_{jk} - (E_{jk}+E_{ik})\delta_{ij}
\end{equation}
whenever $i,j,k$ are distinct. In particular, $\delta_{ik} \in \Z e_i + \Z e_j$. If $n \geq 4$, it follows immediately that $\delta_{ik} \in
\Z e_i$ for all $i \neq k$ and we can easily conclude. In general, it follows that $\delta(1+nE_{ij}) = n \delta_{ij}$. Applying the cocycle relation
to~\eqref{eq.commutator}, it follows that
$$\delta_{ik} = (E_{ij} - E_{ik})\delta_{jk} -
(E_{jk}+nE_{ik})\delta_{ij}$$ whenever $i,j,k$ are distinct and $n \in \Z$. So, $E_{ik} \delta_{ij} = 0$. We conclude that $\delta_{ij} \in
\Z e_i$ for all $i \neq j$. Writing $\delta_{ij} = x_{ij} e_i$ with $x_{ij} \in \Z$,~\eqref{eq.hulpje} reads $x_{ik} = x_{jk}$ for all $i,j,k$
distinct. We have found $y \in \Z^n$ such that $x_{ij} = y_j$ for all $i \neq j$ and then, $\delta(g) = (1-g)y$ for all $g \in \SL(n,\Z)$.

Finally consider a $1$-cocycle $\delta : \GL(n,\Z) \recht \Z^n$, $n \geq 2$. Replacing $\delta$ by a cohomologous $1$-cocycle, we may assume that
$\delta(g) = 0$ whenever $g \in \SL(n,\Z)$. If now $\det T = -1$ and if $g_0 \in \SL(n,\Z)$ is an arbitrary element for which $1$ is not an
eigenvalue, we get
$$0 = \delta(Tg_0 T^{-1}) = (1-Tg_0 T^{-1})\delta(T) \; .$$
So, $\delta(T) = 0$ as well.

So far, we have proven the statements about $\Out(G)$ and $\Out(G^+)$. If $n \geq 3$, $\SL(n,\Z)$ equals its commutator subgroup and so, any
character on it vanishes. This yields the results for $\Char(G)$ and $\Char(G^+)$.

In order to see that $G$ and $G^+$ have infinite conjugacy classes, we look at the conjugacy classes of $(x_0,g_0)$ for $g_0 \neq 1$ and of $(x_0,1)$
for $x_0 \neq 0$. In the first case, all $(x_0 + (1-g_0)y,g_0)$, $y \in \Z^n$ belong to this conjugacy class. In the second case, we of course have
$(G^+ x_0,1)$ in the conjugacy class.
\end{proof}

Denote $J_n = \Bigl(\begin{smallmatrix} 0 & I_n \\ -I_n & 0 \end{smallmatrix} \Bigr) \in \GL(2n,\Z)$. Define the symplectic group $\Sp(2n,\Z)$ as the
elements $A \in \GL(2n,\Z)$ satisfying $A^t J_n A = J_n$. Then, $\Sp(2n,\Z) \subset \SL(2n,\Z)$ and, for $n \geq 2$, $\Sp(2n,\Z)$ is generated by
\begin{equation}\label{eq.generation}
\begin{pmatrix} A & 0 \\ 0 & (A^t)^{-1} \end{pmatrix} \; , \quad \begin{pmatrix} I_n & B \\ 0 & I_n \end{pmatrix} \; , \quad \begin{pmatrix} I_n & 0
\\ C & 0 \end{pmatrix} \quad\text{for}\quad A \in \GL(n,\Z), \quad B,C \in \M_n(\Z)^{\text{sym}}
\end{equation}
where $\M_n(\Z)^{\text{sym}}$ denotes the (additive) group of $n$ by
$n$ integer matrices $B$ satisfying $B = B^t$. See \cite{HOM} for details.

\begin{proposition} \label{prop.outer-two}
Let $n \geq 2$ and write $G = \Z^{2n} \rtimes \Sp(2n,\Z)$. Then, $G$ is an ICC group. The outer automorphism group $\Out(G)$ has
order $2$ and is given by conjugation with an element $A \in
\GL(2n,\Z)$ satisfying $A^t J_n A = - J_n$. If $n \geq 3$, $\Char G$
is trivial. If $n = 2$, $\Char G$ has order $2$.
\end{proposition}
\begin{proof}
It is clear that $G$ is an ICC group. Moreover, $\Sp(2n,\Z)$ is generated by
symplectic transvections, which makes it easy to compute the
commutator subgroup (see \cite{HOM}).

Let now $\al \in \Aut(G)$. As in the proof of Proposition \ref{prop.outer}, it follows that $\al(\Z^{2n}) = \Z^{2n}$, using this time that a
non-zero, $\Sp(2n,\Z)$-invariant subgroup of $\Z^{2n}$ has finite index. But then $\al(x,1) = (X x,1)$ for all $x \in \Z^{2n}$ and some $X \in \GL(2n,\Z)$.
It follows that $X$ normalizes $\Sp(2n,\Z)$, which readily implies that $X^t J_n X = \pm J_n$. For the rest of the proof, we may assume that $\al(x,1)
= (x,1)$ for all $x \in \Z^{2n}$. Hence, $\al(x,g) = (x + \delta(g),g)$, where $\delta : \Sp(2n,\Z) \recht \Z^{2n}$ is a $1$-cocycle (i.e.\ satisfying
the relation $\delta(gh) = \delta(g) + g \delta(h)$ for all $g,h \in \Sp(2n,\Z)$). Since
\begin{equation}\label{eq.special}
g:= \begin{pmatrix} A & 0 \\ 0 & (A^{-1})^t
\end{pmatrix} \in \Sp(2n,\Z) \quad\text{for all} \quad A \in \GL(n,\Z)
\end{equation}
and since all $1$-cocycles $\GL(n,\Z) \recht \Z^n$ are trivial as shown in the proof of Proposition \ref{prop.outer}, we may assume that $\delta(g) =
0$ for all $g$ of the form~\eqref{eq.special}. Define $\delta_i : \M_n(\Z)^{\text{sym}} \recht \Z^n$ by the formula
$$\delta\begin{pmatrix} I_n & B \\ 0 & I_n \end{pmatrix} = \begin{pmatrix} \delta_1(B) \\ \delta_2(B) \end{pmatrix} \; .$$
The $1$-cocycle relation implies that $\delta_2$ is a group homomorphism and $\delta_1(B + C)= \delta_1(B) + \delta_1(C) + B \delta_2(C)$. So,
$B\delta_2(C) = C \delta_2(B)$ for all $B,C \in
\M_n(\Z)^{\text{sym}}$. It follows that $\delta_2(B) = 0$ for all
$B$. But then, $\delta_1$ is a group homomorphism. Conjugating with an element $g$ as in~\eqref{eq.special} and using the $1$-cocycle relation, it follows that $\delta_1(A B A^t) = A
\delta_1(B)$. By induction on $n$, it is checked that such a $\delta_1$ is identically zero. Analogously, $\delta$ is zero on elements of the form
$\Bigl(\begin{smallmatrix} I_n & 0 \\ B & I_n \end{smallmatrix}\Bigr)$. Since the elements in~\eqref{eq.generation} generate $\Sp(2n,\Z)$, it follows
that $\delta = 0$. So, we are done.
\end{proof}

\subsection{Groups satisfying conditions $\cB$ and its stronger versions} \label{subsec.examples-groups}

Recall that a group is called weakly rigid if it admits an infinite
normal subgroup with the relative property~(T).

\begin{examples} \label{ex.conditionsCandD}
\begin{enumerate}
\item Let $G$ be a weakly rigid group and $G_0$ a subgroup
  with the following relative ICC property
$$\{h g h^{-1} \mid h \in G_0 \} \;\;\text{is infinite for all}\;\; g
\in G, g \neq e \; .$$ Then, the diagonal inclusion $G_0 \subset G \times G_0$ satisfies condition $\cB$. The associated generalized Bernoulli action
is defined by the (left-right) double shift $(G \times G_0) \actson G$. Of course, we can take $G_0=G$ whenever $G$ is a weakly rigid ICC group, like
$\Z^n \rtimes \SL(n,\Z)$ ($n \geq 2$), $\PSL(n,\Z)$ ($n \geq 3$) or the direct product of any of these with an arbitrary ICC group.

\item \label{item.two} Let an infinite group $\Gamma$ act on an infinite group $H$ by automorphisms. Assume that $\Gamma \recht \Aut(H)$ is injective and that $H
\subset H \rtimes \Gamma$ has the relative property~(T). We set $G:=H \rtimes \Gamma$ and consider $G_0 = \Gamma \subset G$. The associated
generalized Bernoulli action is defined by the action $G \actson H$ given by $(a,g) \cdot b = a \al_g(b)$ (which satisfies the freeness condition
\ref{conv.free}). If for every $a \in H$, $a \neq e$, $\Stab a$ is of infinite index in $\Gamma$, the group $G$ is ICC and the inclusion $\Gamma
\subset G$ satisfies condition $\cB$. Examples include $\Z^n  \rtimes \SL(n,\Z)$ and $\Z^n \rtimes \GL(n,\Z)$ for all $n \geq 2$. Many more examples
of this type are provided by \cite{TF} and \cite{AV}.

\item \label{item.crucial} One can slightly modify the construction of
  the previous item. Let an infinite group $\Gamma$ act on an infinite
  group $H$ by automorphisms and assume that $\Gamma \recht \Aut(H)$
  is injective. Let $G = H \rtimes \Gamma$ and suppose that $G$ has
  the property~(T). Let $G_0 \subset \Gamma$ such that $gG_0g^{-1}
  \cap G_0$ is finite for all $g \in \Gamma \setmin G_0$.
Suppose that $G_0$ has the Haagerup property and that $\Stab a \cap G_0$ is finite for all $a \in H$, $a \neq e$. Then, $(G,G_0)$ satisfies condition $\cD$. Taking any of the $G_0 \subset \GL(n,\Z)$ as in \ref{ex.malnormal} below, we get explicit examples
of inclusions $G_0 \subset \Z^n \rtimes \GL(n,\Z)$, $n \geq 3$, satisfying condition $\cD$.

\item Taking any of the $G_0$ in \ref{ex.malnormal} below, we consider $G_0 \subset \PSL(n,\Z)$, $n \geq 3$, providing other examples of inclusions satisfying
condition $\cD$.
\end{enumerate}
\end{examples}

We give a construction procedure of subgroups $G_0$ of the groups
$\SL(n,\Z)$ and $\Z^n \rtimes \SL(n,\Z)$ with the property that $gG_0g^{-1} \cap
G_0$ is finite whenever $g \not\in G_0$. Concrete examples are given in \ref{ex.malnormal} below.

Let $\Q \subset K$ be a field extension of degree $n$. Denote by $\cO_K$ the ring of algebraic integers in $K$ and by $\cO_K^\ast$ the group of units of
$\cO_K$. We construct an injective homomorphism $\pi : \cO_K^\ast \rtimes \Gal(K/\Q) \recht \GL(n,\Z)$. Writing $G_0 = \Im \pi$, we prove that
$gG_0g^{-1} \cap G_0$ is finite whenever $g \in \GL(n,\Z)$ and $g \not\in G_0$, provided that there are \emph{no intermediate fields} strictly
between $\Q$ and $K$.

Let $\om_1=1,\om_2,\ldots,\om_{n}$ be a $\Z$-basis of $\cO_K$ and write $\om = (\om_1,\ldots,\om_{n})^t \in \C^n$. Whenever $\si : K \recht \C$ is an
embedding, we write $\om^\si := (\si(\om_1),\ldots,\si(\om_{n}))^t$. Writing the transition matrix between two bases of $\cO_K$, we get
$$\pi : \cO_K^\ast \rtimes \Gal(K/\Q) \recht \GL(n,\Z) : \pi(u,\si)\om = u
\om^\si \quad\text{for}\quad u \in \cO_K^\ast, \si \in \Gal(K/\Q) \; .
$$
Suppose now that there are \emph{no intermediate fields} strictly between $\Q$ and $K$. Write $G_0 = \Im \pi$. In order to show that $gG_0g^{-1} \cap
G_0$ is finite whenever $g \in \GL(n,\Z)$ and $g \not\in G_0$, it suffices to show the following: if $u,w \in \cO_K^\ast \setmin \{\pm 1\}$, $g \in
\GL(n,\Z)$ and $g \pi(u) g^{-1} = \pi(w)$, then $g \in G_0$. If $\si_1,\ldots,\si_n$ are the $n$ embeddings of $K$ into $\C$, write
$\om^{(i)}:=\om^{\si_i}$ and $u^{(i)}:= \si_i(u)$ for $u \in K$. Whenever $u \in \cO_K^\ast \setmin \{\pm 1\}$, $u$ generates $K$ as a field, by
our assumption. It follows that the $u^{(i)}$ are distinct. Moreover, $\pi(u) \om^{(i)} = u^{(i)} \om^{(i)}$ and we conclude that
$\{\om^{(1)},\ldots,\om^{(n)}\}$ is a basis of eigenvectors of $\pi(u)$, with distinct eigenvalues $u^{(1)},\ldots,u^{(n)}$. If now $u,w \in
\cO_K^\ast \setmin \{\pm 1\}$, $g \in \GL(n,\Z)$ and $g \pi(u) g^{-1} = \pi(w)$, it follows that $g \om = \gamma \om^{\si}$ for $\gamma \in \C$ and
$\si : K \recht \C$ an embedding. It follows that $\gamma \in \cO_K^\ast$ and $\si \in \Gal(K/\Q)$. So, $g \in G_0$.

\begin{examples} \label{ex.malnormal}
The methods of computational algebraic number theory
  (and their implementation in PARI/GP) allow to obtain as many
  concrete examples as you like.
\begin{enumerate}
\item Let $G = \GL(2,\Z)$ and $G_0 = \{ \pm \bigl(\begin{smallmatrix}
    5 & -1 \\ 1 & 0 \end{smallmatrix}\bigr)^n \mid n \in \Z \}$. Then,
    $gG_0g^{-1} \cap G_0 = \{\pm 1\}$ for all $g \in G \setmin G_0$.
\item \label{item.three} Let $G = \SL(3,\Z)$. Using the extension by joining a root of
  the polynomial $x^3+x+1$, we define
$$A=\begin{pmatrix} 0 & -1 & -1 \\ 1 & 0 & 0 \\ 0 & 1 & 0 \end{pmatrix}
\quad\text{and}\quad G_0 = A^\Z \; .$$ Then, $gG_0g^{-1} \cap G_0 = \{1\}$ whenever $g \in G \setmin G_0$.
\item \label{item.four} Let $G = \GL(4,\Z)$. Using the extension by joining a root of
  the polynomial $x^4 + x + 1$, we define
\begin{equation}
A=\begin{pmatrix} 0 & 0 & -1 & -1 \\ 1 & 0 & 0 & 0 \\ 0 & 1 & 0 & 0
\\ 0 & 0 & 1 & 0\end{pmatrix}
\quad\text{and}\quad G_0 = \{\pm A^n \mid n \in \Z\} \; . \label{eq.thisoneworks}
\end{equation}
Then, $gG_0g^{-1} \cap G_0 = \{\pm 1\}$ whenever $g \in G \setmin G_0$.
\end{enumerate}
\end{examples}

\begin{example} \label{ex.no-anti-aut}
Let $n \geq 3$. The methods of \ref{ex.malnormal} allow to write a
subgroup $H \subset \GL(n,\Z)$ with the following properties (below we
give a concrete example).
\begin{itemize}
\item $H$ has a finite index normal subgroup $H_0$.
\item $H_0$ is an abelian group of symmetric matrices.
\item Two elements of $H_0$ that share an eigenvalue are conjugate in
  $H$.
\item If $g \in H_0$ and $g \neq 1$, then $g$ has $n$ distinct
  eigenvalues and none of them is an eigenvalue of $g^{-1}$.
\item If $g \in H_0$ and $g \neq 1$, the centralizer of $g$ in
  $\GL(n,\Z)$ is contained in $H$.
\end{itemize}
It follows that whenever $g \in \GL(n,\Z)$ and $gHg^{-1} \cap H$
infinite, then $g \in H$. Indeed, in such a case, also $gH_0g^{-1}
\cap H_0$ is infinite and we take $h \in H_0 \setmin \{e\}$ such
that $g h g^{-1} \in H_0$. Our assumptions imply that $g h g^{-1}$ is
conjugate to $h$ in $H$. Multiplying $g$ by an element of $H$, we may
assume that $g h g^{-1} = h$. But then, $g \in H$.

Denote by $\al$ the automorphism of $\GL(n,\Z)$ defined by $\al(g) =
(g^{-1})^t$. It also follows that $\al(H) = H$. Indeed, by assumption,
$\al(H_0) = H_0$ since $H_0$ consists of symmetric matrices. If then
$g \in H$, $\al(g)$ normalizes $H_0$, implying that $\al(g) \in H$.

Define $\Gamma = \Sp(2n,\Z)$. Write $J$ instead of $J_n$. Denote $\pi : \GL(n,\Z) \recht \Gamma : \pi(g) = \bigl(
\begin{smallmatrix} g & 0 \\ 0 & \al(g)
\end{smallmatrix}\bigr)$. Define $\Gamma_0$ as the subgroup of
$\Gamma$ generated by $J$ and $\pi(H)$. We claim that whenever $g \in
\Gamma$ and $g \Gamma_0 g^{-1} \cap \Gamma_0$ infinite, then $g \in \Gamma_0$.

Let $g \in
\Gamma$ and $g \Gamma_0 g^{-1} \cap \Gamma_0$ infinite. Note that $J
\pi(g) J^{-1} = \pi(\al(g))$. So, $\pi(H_0)$ is a finite index normal
subgroup of $\Gamma_0$. Hence, $g \pi(H_0) g^{-1} \cap \pi(H_0)$ is
infinite. Take $a,b \in H_0 \setmin \{e\}$ such that $g \pi(a)
g^{-1}  = \pi(b)$. If $\spect(a)$ denotes the set of eigenvalues of $a$,
we have $\spect(\pi(a)) = \spect(a) \cup \spect(a)^{-1}$. Hence, either $a,b$,
either $a,b^{-1}$ share an eigenvalue. In the first case, we replace
$g$ by $\pi(h)g$ for a suitable element $h \in H$ and in the second
case by $\pi(h)Jg$. In both cases, we may assume that $g$ commutes
with $\pi(a)$. Since $\pi(a)$ has $2n$ distinct eigenvalues, it
follows that $g$ has the form $\bigl(\begin{smallmatrix} A & 0 \\ 0 &
  B \end{smallmatrix}\bigr)$. Since $g \in \Sp(2n,\Z)$, it follows
that $g = \pi(h)$ for some $h \in \GL(n,\Z)$. But then, $h$ and $a$
commute, implying that $h \in H$ and hence, $g \in \Gamma_0$.

If we write $G = \Z^{2n} \rtimes \Sp(2n,\Z)$ with $G_0 = \Gamma_0$
viewed as a subgroup of $G$, it follows that the pair $(G,G_0)$
satisfies condition $\cD$. It can indeed by easily checked that $G$
has property~(T): if $n \geq 2$, $\Sp(2n,\Z)$ has property~(T) and
$\Z^{2n} \subset \Z^{2n} \rtimes \Sp(2n,\Z)$ has the relative property
(T). The latter follows because it contains $(\Z^n
\oplus \Z^n) \rtimes \SL(n,\Z)$ with $\SL(n,\Z)$ acting diagonally.

We finally give a concrete example of $H \subset \GL(3,\Z)$ satisfying
the above requirements. Let $x$ be a root of the equation $x^3-3x+1=0$
and $K$ the subfield of $\R$ generated by $x$. Then, $1$, $x$ and
$x^2$ form a $\Z$-basis of $\cO_K$. The group $\cO_K^*$ is
isomorphic with $\Z^2 \oplus \Z/2\Z$ and generated by $x,x-1$ and
$-1$. Then, $x,-1+x+x^2,1$ is as well a $\Z$-basis of $\cO_K$,
yielding the homomomorphism $\pi : \cO_K^* \recht \GL(3,\Z)$ given by
\begin{equation}\label{eq.onzematrices}
\pi(x) = \begin{pmatrix} -1 & 1 & 1 \\ 1 & 1 & 0 \\ 1 & 0 & 0
\end{pmatrix} \; , \quad \pi(x-1) =  \begin{pmatrix} -2 & 1 & 1 \\ 1 & 0 & 0 \\ 1 & 0 & -1
\end{pmatrix} \; , \quad \pi(-1) = -I \; .
\end{equation}
It follows that $\pi(\cO_K^*)$ consists of symmetric matrices and we
define $H := \pi(\cO_K^* \rtimes \Gal(K/\Q))$. For completeness, we
mention that $\Gal(K/\Q) \cong \Z/3\Z$, yielding that $H$ is generated
by the matrices in~\eqref{eq.onzematrices} together with the matrix of
order $3$ given by $\Bigl(\begin{smallmatrix} -1 & 1 & -1 \\ -1 & 0 &
  1 \\ 0 & 0 & 1 \end{smallmatrix}\Bigr)$.
\end{example}

We are now ready to prove the following.

\begin{proposition} \label{prop.everygroup}
Let $Q$ be a group of finite presentation. There exists $(G,G_0)$ in the family $\cF_3$ introduced in Proposition \ref{prop.families}, satisfying $\Char(G) =
\{1\}$ and
$$\frac{\Aut(G_0 \subset G)}{\Ad G_0} \cong Q \; .$$
\end{proposition}

\begin{lemma} \label{lem.computeOut}
Let $H$ be a non-amenable group, with $\Char(H) = \{1\}$. Suppose that
the centralizer $C_H(h)$ is amenable for all $h \in H \setmin \{e\}$. Let the finite group
$\Gamma$ act faithfully on the finite set $I$ (i.e.\ $\Gamma \recht
\operatorname{Perm}(I)$ is injective). We set
$$G = K \times (H^I \rtimes \Gamma) \quad\text{and}\quad G_0 = K_0
\times \Delta(H) \times \Gamma \; ,$$
where $K$ is an ICC group with amenable subgroup $K_0$ and
where $\Delta : H \recht H^I$ is the diagonal embedding of $H$ into
$H^I := \oplus_I H$. Then,
\begin{itemize}
\item The natural map yields an isomorphism $\dis \frac{\Aut(G_0 \subset G)}{\Ad G_0} \cong \frac{\Aut(K_0
    \subset K)}{\Ad K_0} \times \Out(H) \times \frac{\Aut^*(\Gamma
    \actson I)}{\Gamma}$.
\item $G$ is ICC.
\item $\Char(G) \cong \Char(K) \times \Char(\Gamma)$.
\end{itemize}
\end{lemma}
\begin{proof}
Only the first statement is non-trivial.
Let $\al \in \Aut(G_0 \subset G)$. Consider $p_{K_0 \times \Gamma} :
G_0 \recht K_0 \times \Gamma$. It follows that the kernel of $(p_{K_0
  \times \Gamma} \circ \al)|_{\Delta(H)}$ is non-amenable yielding a
non-amenable subgroup $H_0 \subset H$ such that $\al(\Delta(H_0))
\subset \Delta(H)$. By our assumption on the centralizers in $H$, it
follows that $\al(K \times \Gamma) \subset K \times \Gamma$ and by
symmetry, the equality holds. Since the ICC group $K$ does not have
finite normal subgroups, $\al(\Gamma) = \Gamma$. We can take as well a
finite index subgroup $L \subset K$ such that $\al(L) \subset
K$. Since $K$ is an ICC group, the centralizer of $L$ in $K$ is
trivial and it follows that $\al(H^I) \subset H^I
\rtimes \Gamma$. Because $\al(\Gamma) = \Gamma$,
$p_\Gamma(\al(\Delta(H)))$ belongs to the center of $\Gamma$.
Since $\Char H = \{1\}$, the homomorphism $h \in H \mapsto
p_\Gamma(\al(\Delta(h)))$ is trivial.
So, $\al(\Delta(H)) =
\Delta(H)$. Whenever $g \in H$ and $i \in I$, we write $g_i := (e,\ldots,g,\ldots,e)$ where $g$ appears in position $i$. Since $p_\Gamma(\al(\Delta(h))) = e$ for
all $h \in H$, it follows that
$$p_\Gamma(\al((ghg^{-1}h^{-1})_i)) = p_\Gamma(\al(g_i \Delta(h) g_i^{-1} \Delta(h)^{-1})) = e \; .$$
Since $H$ equals its commutator subgroup, we get $\al(H^I) = H^I$. We
already got $\al(\Gamma) = \Gamma$ and it follows that $\al(K) = K$.

Using amenability of the centralizers $C_H(h)$ for $h \in H \setmin \{e\}$, it
is easy to check that every $\al \in \Aut(\Delta(H) \subset H^I)$ is of the form $\al(h_i) = (\al_0(h))_{\sigma(i)}$ for all $h \in H$, $i \in I$,
where $\al_0 \in \Aut(H)$ and $\sigma : I \recht I$ a permutation. So,
we have shown that $\al \in \Aut(K_0 \subset K) \times \Aut(H) \times \Aut^*(\Gamma \actson I)$.
\end{proof}

\begin{lemma} \label{lem.finitegroup}
Denote by $F_2$ the field of order $2$. For $n \geq 3$, set $\Gamma = \GL(n,F_2)$. Then, $\Aut^*(\Gamma \actson F_2^n) = \Gamma$ and $\Char \Gamma =
\{1\}$.
\end{lemma}
\begin{proof}
The group $\Gamma$ has no characters and one outer automorphism given
by $\al(g) = (g^{-1})^t$ (see \cite{Dieud}). We claim that
$\Aut^\al(\Gamma \actson F_2^n)$ is empty. Let $e_1$ be the first
basis vector of $F_2^n$. If $\si \in \Aut^\al(\Gamma \actson F_2^n)$,
it follows that $\si(e_1)$ is fixed by all the matrices in $\al(\Stab
e_1)$. Such an element does not exist in $F_2^n$.

It remains to show that any permutation $\si$ of $F_2^n$ commuting
with the action of $\Gamma$, is the identity map. But this is obvious
since for $e \in F_2^n$, $\Stab e$ fixes only $e$ itself.
\end{proof}

\begin{lemma} \label{lem.TandTriv}
Set $K = \Z^4 \rtimes \SL(4,\Z)$ and $K_0 \subset \SL(4,\Z) \subset K$ consisting of the elements $\pm A^n$, $n \in \Z$, where $A$ is given by
\ref{ex.malnormal}.\ref{item.four}. Then, $K$ is an ICC group,
$\Char(K) = \{1\}$ and $\Aut(K_0 \subset K) = \Ad K_0$. Moreover, $K_0
\subset K$ is almost malnormal.
\end{lemma}
\begin{proof}
By \ref{ex.conditionsCandD}.\ref{item.crucial}, $K_0$ is almost
malnormal in $K$.
By Proposition \ref{prop.outer}, the only outer automorphism of $K$ is
given by conjugation by an element $T \in \GL(4,\Z)$ with determinant
$-1$. From \ref{ex.malnormal}.\ref{item.four}, it follows that for
such a $T$, we have $TK_0T^{-1} \neq K_0$. We conclude that $\Aut(K_0
\subset K) = \Ad K_0$. Finally, $\Char K = \{1\}$ by Proposition \ref{prop.outer}.
\end{proof}

We are then ready to prove Proposition \ref{prop.everygroup}.

\begin{proof}[Proof of Proposition \ref{prop.everygroup}]
Using \cite{Bum-Wise}, we take a finitely generated group $H$ with the following properties: $\Out(H) \cong Q$ and $H$ is a non-elementary
subgroup of a finitely presented $C'(\frac{1}{6})$ small cancellation group. Slightly adapting the construction of \cite{Bum-Wise}, we may assume that $\Char(H) =
\{1\}$. Small cancellation implies that $C_H(h)$ is amenable (even
cyclic) for every $h \in H \setmin \{e\}$. Since a finitely presented $C'(\frac{1}{6})$
small cancellation group is word hyperbolic, the subgroup $H$ is in
Ozawa's class $\cC$ (see \cite{OzSolid,OzKurosh}) and there exist $h_1,h_2$ such that $C_H(h_1,h_2)$ is trivial.

Let $K_0 \subset K$ be as in Lemma \ref{lem.TandTriv} and $\Gamma \actson I$ as in Lemma \ref{lem.finitegroup}. Set
$$G = K \times (H^I \rtimes \Gamma) \quad\text{and}\quad G_0 = K_0 \times \Delta(H) \times \Gamma \; ,$$
where, as above, $\Delta : H \recht H^I$ is the diagonal embedding. From Lemma \ref{lem.computeOut}, it follows that $\frac{\Aut(G_0 \subset G)}{\Ad
G_0} \cong Q$. It is readily checked that for $g \in G \setmin G_0$, we have either $gG_0g^{-1} \cap G_0$ amenable (in the case where $g \not\in K
\times \Delta(H) \times \Gamma$), or $gG_0g^{-1} \cap G_0 \cong \Lambda \times H$ for some finite group $\Lambda$ (in the case where $g \in (K-K_0) \times
\Delta(H) \times \Gamma$). In all cases, $gG_0g^{-1} \cap G_0$ belongs to the class $\cC$ of Ozawa.

Using $h_1,h_2$ and some element in $K \setmin K_0$, it is easy to write elements $g_1,\ldots,g_n \in G$ such that $\bigcap_{i=1}^n g_i G_0
g_i^{-1}$ is finite.
\end{proof}

It is easier to give examples of $(G,G_0)$ satisfying condition $\cB$ (Def.\ \ref{def.conditionB}) such that $\frac{\Aut(G_0 \subset G)}{\Ad G_0}
\cong Q$ for an arbitrary countable group $Q$.

\begin{proposition} \label{prop.everygroup-bis}
Let $Q$ be any countable group. There exists $(G,G_0)$ satisfying condition $\cB$ (Def.\ \ref{def.conditionB}) such that $\frac{\Aut(G_0 \subset
G)}{\Ad G_0} \cong Q$. Moreover, $(G,G_0)$ can be chosen in such a way that every character on $G$ that is trivial on $G_0$, is trivial on the whole
of $G$.
\end{proposition}
\begin{proof}
Let $K_0 \subset K$ be as in Lemma \ref{lem.TandTriv}. Let $H$ be a
non-amenable ICC group without characters and with $C_H(h)$ amenable
for all $h \in H-\{e\}$. Set
$$G= K \times bigl((H \times H) \rtimes \frac{\Z}{2\Z}\bigr) \quad\text{and}\quad G_0 = K_0 \times \Delta(H) \times
\frac{\Z}{2\Z}\; ,$$
where $\Z/2\Z$ acts on $H \times H$ by the flip automorphism and where $\Delta : H \recht H \times H$ is the diagonal embedding.
Then, $(G,G_0)$ satisfies condition $\cB$ and
$$\frac{\Aut(G_0 \subset G)}{\Ad G_0} \cong \Out(H)\; .$$
The latter can be any countable group, by the results of
\cite{Bum-Wise}. Moreover, any character on $G$ that is trivial on
$G_0$, is trivial on the whole of $G$.
\end{proof}

\subsection{Outer automorphisms of equivalence relations} \label{subsec.outer-equivalence}

Let $n \geq 2$ and consider $G:=\Z^n \rtimes \GL(n,\Z)$ acting on $\Z^n$ by $(x,g) \cdot y = x+g y$, which is $G \actson G/G_0$ for $G_0 =
\GL(n,\Z)$. By \ref{ex.conditionsCandD}.\ref{item.two}, the inclusion $G_0 \subset G$ satisfies condition $\cB$. Let $\cR$ be the type II$_1$
equivalence relation given by the orbits of
$$G \actson \Bigl( \prod_{\Z^n} (X_0,\mu_0) \Bigr)^K$$
whenever $K \actson (X_0,\mu_0)$ is a faithful action of a compact group $K$ on the base $(X_0,\mu_0)$. By Proposition \ref{prop.outer}, $\Aut(G_0
\subset G) = \Ad G_0$. Also $\Hom(G/G_0 \recht \cZ(K))$ is trivial, because every character of $G$ that is trivial on $G_0$ is itself trivial. So,
Theorem \ref{thm.classification} yields
$$\Out(\cR) \cong \frac{\Aut^*(K \actson (X_0,\mu_0))}{K} \; .$$
As observed in the introduction, we can take $K=\{e\}$ and $(X_0,\mu_0)$ atomic with distinct weights, producing a continuous family of non
stably isomorphic equivalence relations with trivial outer automorphism group.

But, another interesting case is to take the action $K \actson K$ for which $$\frac{\Aut^*(K \actson K)}{K} \cong \Aut K$$ in the obvious way. We can
obtain continuous families of non stably isomorphic equivalence relations as follows. Take $(X_0,\mu_0)$ to be a disjoint union of copies of
$K$ with different mass for $\mu_0$. One still has $\frac{\Aut^*(K \actson K)}{K} \cong \Aut K$. Applying Proposition \ref{prop.everygroup-bis},
we arrive at the following.

\begin{theorem} \label{thm.lot-of-out-equiv}
Let $Q$ be any countable group and $K$ any second countable compact group. There exists a continuous family (explicitly constructed above) of type
II$_1$ equivalence relations $\cR$ that are non stably isomorphic and satisfy
$$\Out(\cR) \cong Q \times \Aut(K) \; .$$
\end{theorem}

\subsection{Outer automorphisms of II$_1$ factors} \label{subsec.outer-vN}

The discussion in Subsection \ref{subsec.examples-groups}, immediately yields the following continuous family of II$_1$ factors $M$ with $\Out(M)$
trivial and with trivial fundamental group. Set $G = \Z^4 \rtimes \SL(4,\Z)$ and $G_0 \subset \SL(4,\Z) \subset G$ consisting of the elements $\pm
A^n$, $n \in \Z$, where $A$ is given by \ref{ex.malnormal}.\ref{item.four}. The outer automorphism group of the II$_1$ factor $L^\infty(\prod_{G/G_0}
(X_0,\mu_0)) \rtimes G$ is given by $\Aut(X_0,\mu_0)$ and the II$_1$ factor remembers $(X_0,\mu_0)$.

Applying Proposition \ref{prop.everygroup}, we get the following.

\begin{theorem} \label{thm.lot-of-out}
Let $Q$ be any finitely presented group. There exists a continuous family (explicitly constructed above) of II$_1$ factors $M$ satisfying $\Out(M) \cong Q$,
$\cF(M) = \{1\}$ and $\Inn M$ open in $\Aut M$.
\end{theorem}

Obviously, if $M$ is a II$_1$ factor with separable predual and with $\Inn M$ open in $\Aut M$, we have $\Out M$ a countable group and the theorem
shows that any finitely presented group can arise in this
way. Proposition \ref{prop.everygroup-bis} suggests that one can
actually obtain any countable group.
If we only
assume that $\Inn M$ is closed in $\Aut M$ (such an $M$ is said to be \emph{full}), the group $\Out M$ is a Polish group. We
do not know which Polish groups can arise as $\Out M$ for a full
II$_1$ factor $M$. By
\cite{IPP}, all abelian compact second countable groups do arise in
this way.

\begin{proposition} \label{prop.millenium}
Suppose that $G_0 \subset G$ belongs to one of the families in \ref{prop.families}. Define the wreath product group
$$H := \Bigl(\bigoplus_{G/G_0} \frac{\Z}{2\Z} \Bigr) \rtimes G \; .$$
Then, the natural map yields an isomorphism
$$\Out(\cL(H)) \cong \Char(H) \rtimes \Out(H) \; .$$
\end{proposition}
\begin{proof}
To every $(\om,\delta) \in \Char(H) \rtimes \Aut(H)$, we associate the automorphism of $\cL(H)$ defined by $\theta_{\om,\delta}(u_g) = \om(g)
u_{\delta(g)}$. Since $H$ is an ICC group, we get an injective homomorphism $\Char(H) \rtimes \Out(H) \recht \Out(\cL(H))$. Observe that $\cL(H)$ is
naturally identified with the crossed product of the $(G \actson G/G_0)$-Bernoulli action with base $(\cL(\Z/2\Z),\tau)$. Since the only non-trivial
automorphism of $\cL(\Z/2\Z)$ can be written as $\theta_\om(u_g) = \om(g) u_g$, where $\om$ is the non-trivial character of $\Z/2\Z$, it follows from
Theorem \ref{thm.vNstrong} that, up to inner automorphisms, every automorphism of $\cL(H)$ is of the form $\theta_{\om,\delta}$.
\end{proof}

\begin{theorem} \label{thm.without-anti}
Let $G = \Z^6 \rtimes \Sp(6,\Z)$ and $J = \bigl(\begin{smallmatrix} 0
  & I \\ -I & 0 \end{smallmatrix}\bigr) \in \GL(6,\Z)$. Let $G_0
  \subset G$ be the subgroup
defined in \ref{ex.no-anti-aut}. For every $\al \in \, ]0,\pi/2[$,
define the $2$-cocycle $\Om_\al \in Z^2(G,S^1)$ by
$$\Om_\al(x,y) = \exp(i \al \, x^t J y) \quad\text{for}\quad x,y \in
\Z^6\; , \quad \Om_\al((x,g),(y,h)) := \Om_\al(x,g\cdot y) \; .$$
For every standard probability space $(X_0,\mu_0)$, consider the $(G
\actson G/G_0)$-Bernoulli action with base $(X_0,\mu_0)$ and write
$$M_\al(X_0,\mu_0) := L^\infty(X,\mu) \twist{\Om_\al} G \; .$$
Then, $M_\al(X_0,\mu_0)$ has trivial fundamental group, has no
anti-automorphism and has outer automorphism group given by
$\Aut(X_0,\mu_0)$. Moreover, $M_\al(X_0,\mu_0)$ and
$M_\beta(Y_0,\eta_0)$ are stably isomorphic if and only if $\al=\beta$
and $(X_0,\mu_0) \cong (Y_0,\eta_0)$.
\end{theorem}
\begin{proof}
The theorem follows from Theorem \ref{thm.vNstrong-cocycle}, using
Example \ref{ex.no-anti-aut} and Proposition \ref{prop.outer-two} and
observing that the formula
$$\Om_\al(x,y) = \exp(i \al \, x^t J y)$$
defines, for every $\al \in \R$, a $2$-cocycle $\Om_\al \in
Z^2(\Z^n,S^1)$. Moreover $\Om_\al$ is a coboundary if and only if
$\al \in \Z\pi$. Our choice of $\al \in \, ]0,\pi/2[$ guarantees that
$\Om_\al$ and $\overline{\Om_\al}$ define different elements of
$H^2(G,S^1)$. Then, $\Aut_{\Om_\al}(G) = \Ad G$ as well.
\end{proof}


\begin{thebibliography}{AA}
\bibitem{Borel} {\sc A. Borel},
On the automorphisms of certain subgroups of semi-simple Lie
groups. In {\it Algebraic Geometry, Tata Inst. Fund. Res., Bombay,
  1968}, Oxford Univ. Press, London, pp. 43--73.

\bibitem{Bum-Wise} {\sc I. Bumagin \& D.T. Wise}, Every group is an outer automorphism group of a finitely generated group. {\it J. Pure Appl. Algebra} {\bf 200} (2005),
137--147.

\bibitem{C1} {\sc A. Connes}, A factor of type II$_1$ with
  countable fundamental group. {\it J. Operator Theory} {\bf 4}
  (1980), 151--153.

\bibitem{C2} {\sc A. Connes}, Classification of injective factors.
{\it Ann. of Math. (2)} {\bf 104} (1976), 73--115.

\bibitem{C5} {\sc A. Connes}, Classification des facteurs.
In {\it Operator algebras and applications, Part 2 (Kingston, Ont., 1980)},
{\it Proc. Sympos. Pure Math.} {\bf 38}, Amer. Math. Soc., Providence,
1982, p.~43--109.

\bibitem{C10} {\sc A. Connes},
Sur la classification des facteurs de type ${\rm II}$.
{\it C. R. Acad. Sci. Paris S{\'e}r. A} {\bf 281} (1975), 13--15.

\bibitem{CJ} {\sc A. Connes \& V.F.R. Jones}, Property (T) for von
  Neumann algebras. {\it Bull. London Math. Soc.} {\bf 17} (1985), 57--62.

\bibitem{CH} {\sc M. Cowling \& U. Haagerup}, Completely bounded multipliers of the Fourier algebra of a simple Lie group of real rank one.
{\it Invent. Math.} {\bf 96} (1989), 507--549.

\bibitem{Dieud} {\sc J. Dieudonn{\'e}}, La g{\'e}om{\'e}trie des groupes
  classiques. Springer Verlag, Berlin, 1955.

\bibitem{Dye1} {\sc H.A. Dye}, On groups of measure preserving
  transformation. I. {\it Amer. J. Math.} {\bf 81} (1959), 119--159.

\bibitem{Dye2} {\sc H.A. Dye}, On groups of measure preserving
  transformation. II.  {\it Amer. J. Math.} {\bf 85} (1963), 551--576.

\bibitem{FM} {\sc J. Feldman \& C.C. Moore}, Ergodic equivalence relations, cohomology, and von Neumann algebras, II.
{\it Trans. Amer. Math. Soc.} {\bf 234} (1977), 325--359.

\bibitem{TF} {\sc T. Fern{\'o}s}, Relative property~(T) and linear
  groups. {\it Ann. Inst. Fourier}, in press. {\tt math.GR/0411527}

\bibitem{Fur0} {\sc A. Furman}, Outer automorphism groups of some
  ergodic equivalence relations. {\it Comment. Math. Helv.} {\bf 80} (2005), 157--196.

\bibitem{Fur1} {\sc A. Furman}, Gromov's measure equivalence and
  rigidity of higher rank lattices. {\it Ann. of Math. (2)} {\bf
  150} (1999), 1059--1081.

\bibitem{Fur2} {\sc A. Furman}, Orbit equivalence rigidity.
{\it Ann. of Math. (2)} {\bf 150} (1999), 1083--1108.

\bibitem{Gab1} {\sc D. Gaboriau}, Invariants $l^2$ de relations
  d'{\'e}quivalence et de groupes. {\it Publ. Math. Inst. Hautes {\'E}tudes
  Sci.} {\bf 95} (2002), 93--150.

\bibitem{Gab2} {\sc D. Gaboriau}, Co{\^u}t des relations d'{\'e}quivalence et
  des groupes. {\it Invent. Math.} {\bf 139} (2000), 41--98.

\bibitem{Geft1} {\sc S.L. Gefter}, Ergodic equivalence relations without outer automorphisms.
{\it Dokl. Acad. Sci. of Ukraine} {\bf 11} (1993), 8–-10.

\bibitem{Geft2} {\sc S.L. Gefter}, Outer automorphism group of the ergodic equivalence relation
generated by translations of dense subgroup of compact group on its homogeneous space. {\it Publ. RIMS, Kyoto Univ.} {\bf 32} (1996), 517–-538.

\bibitem{GdlH} {\sc E. Ghys \& P. de la Harpe}, editors,
Sur les groupes hyperboliques d'apr{\`e}s Mikhael Gromov. Birkh{\"a}user,
Boston, 1990.

\bibitem{GN} {\sc V. Ya. Golodets \& N.I. Nessonov},
T-property and nonisomorphic full factors of types II and III. {\it J. Funct. Anal.} {\bf 70} (1987), 80--89.

\bibitem{HOM} {\sc A.J. Hahn \& O.T. O'Meara}, The classical groups and $K$-theory. {\it Grundlehren der
Mathematischen Wissenschaften} {\bf 291}, Springer-Verlag, Berlin, 1989.

\bibitem{Hjo} {\sc G. Hjorth}, A converse to Dye's theorem.
{\it Trans. Amer. Math. Soc.} {\bf 357} (2005), 3083--3103.

\bibitem{HR} {\sc L.K. Hua \& I. Reiner},
Automorphisms of the unimodular group.
{\it Trans. Amer. Math. Soc.} {\bf 71} (1951), 331--348.

\bibitem{IPP} {\sc  A. Ioana, J. Peterson \& S. Popa}, Amalgamated
  free products of $w$-rigid factors and calculation of their symmetry
  groups. {\it Preprint.} {\tt math.OA/0505589}

\bibitem{JS} {\sc P. Jolissaint \& Y. Stalder}, Strongly singular and strongly mixing MASA's in finite von Neumann algebras. {\it Preprint.} {\tt
math.OA/0602158}

\bibitem{JoMill} {\sc V.F.R. Jones},
Ten problems. In {\it Mathematics: frontiers and perspectives},
Amer. Math. Soc., Providence, 2000, pp.~79--91.

\bibitem{Marg} {\sc G.A. Margulis},
Discrete subgroups of semisimple Lie groups. Springer Verlag, Berlin, 1991.

\bibitem{MS} {\sc N. Monod \& Y. Shalom}, Orbit equivalence
  rigidity and bounded cohomology. {\it Ann. Math.}, to appear.

\bibitem{OzSolid} {\sc N. Ozawa}, Solid von Neumann algebras.
{\it Acta Math.} {\bf 192} (2004), 111--117.

\bibitem{OzKurosh} {\sc N. Ozawa}, A Kurosh type theorem for type II$_1$ factors.
{\it Int. Math. Res. Not.}, to appear. {\tt math.OA/0401121}

\bibitem{P0} {\sc S. Popa}, Cocycle and orbit equivalence superrigidity
for malleable actions of $w$-rigid groups. {\it Preprint.} {\tt
  math.GR/0512646}

\bibitem{P1} {\sc S. Popa}, Strong rigidity of II$_1$ factors
arising from malleable actions of $w$-rigid groups, Part I. {\it
  Invent. Math.}, to appear.
{\tt math.OA/0305306}

\bibitem{P2} {\sc S. Popa}, Strong rigidity of II$_1$ factors
  arising from malleable actions of $w$-rigid groups, Part II.
  {\it Invent. Math.}, to appear. {\tt math.OA/0407103}

\bibitem{P3} {\sc S. Popa}, Some computations of $1$-cohomology groups
and construction of non orbit equivalent actions. {\it
  J. Inst. Math. Jussieu} {\bf 5} (2006), 309--332.

\bibitem{P4} {\sc S. Popa}, Some rigidity results for non-commutative
Bernoulli shifts. {\it J. Funct. Anal.} {\bf 230} (2006),
273--328.

\bibitem{P5} {\sc S. Popa}, On a class of
type II$_1$ factors with Betti numbers invariants. {\it Ann. of
  Math.} {\bf 163} (2006), 809--899.

\bibitem{P6} {\sc S. Popa}, Orthogonal pairs of $^*$-subalgebras in finite von Neumann algebras.
{\it J. Operator Theory} {\bf 9} (1983), 253--268.

\bibitem{P8} {\sc S. Popa}, Correspondences. {\it INCREST
    Preprint} (1986).

\bibitem{PS} {\sc S. Popa \& R. Sasyk}, On the cohomology of actions of
groups by Bernoulli shifts. {\it Preprint.} {\tt math.OA/0311417}

\bibitem{V} {\sc S. Vaes}, Rigidity results for Bernoulli shifts and their von Neumann algebras
(after Sorin Popa). S{\'e}m. Bourbaki, exp.\ no.\ 961. {\it Ast{\'e}risque},
to appear. {\tt math.OA/0603434}

\bibitem{AV} {\sc A. Valette}, Group pairs with property~(T), from
  arithmetic lattices. {\it Geom. Dedicata} {\bf 112} (2005), 183--196.

\bibitem{Zim1} {\sc R.J. Zimmer}, Strong rigidity for ergodic actions
  of semisimple Lie groups. {\it Ann. of Math. (2)} {\bf 112} (1980), 511--529.

\bibitem{Zim} {\sc R.J. Zimmer}, Ergodic theory and semisimple
    groups, Birkh{\"a}user, 1984.
\end{thebibliography}
\end{document}